\documentclass[11pt]{article}
\usepackage[utf8]{inputenc}
\usepackage[english]{babel}
\usepackage{graphicx}
\usepackage{amsmath}
\usepackage{amssymb}
\usepackage{amsfonts}
\usepackage{amsthm}
\usepackage[left=2.1cm,top=1.5cm,right=2.1cm, bottom=2.2cm,letterpaper]{geometry}
\usepackage{mathrsfs}
\usepackage{caption}
\usepackage{subcaption}
\usepackage{latexsym}
\usepackage{xfrac}
\usepackage{tcolorbox}
\usepackage{enumitem}
\usepackage{float}
\usepackage{hyperref}
\usepackage[all]{hypcap}
 \usepackage{url}
\usepackage[toc]{appendix}
\newtheorem{theorem}{Theorem}[section]
\newtheorem{corollary}{Corollary}[section]
\newtheorem{lemma}{Lemma}[section]
\newtheorem{definition}{Definition}[section]

\newtheorem{remark}{Remark}[section]
\numberwithin{equation}{section}
\numberwithin{figure}{section}

\newcommand{\hlabel}{\phantomsection\label}
\makeatletter
\renewcommand*\env@matrix[1][*\c@MaxMatrixCols c]{%
	\hskip -\arraycolsep
	\let\@ifnextchar\new@ifnextchar
	\array{#1}}
\makeatother

\def\neweq#1{\begin{equation}\label{#1}}
\def\endeq{\end{equation}}

\begin{document}
	
\title{Homogenization of Leray's flux problem for the steady-state Navier-Stokes equations in a multiply-connected planar domain}

\author{Clara PATRIARCA - Gianmarco SPERONE\footnote{Research supported by MUR grant \textit{Dipartimento di Eccellenza 2023-2027}.}}
\date{}
\maketitle
\vspace*{-6mm}
\begin{abstract}
	\noindent
The steady motion of a viscous incompressible fluid in a multiply-connected, planar, bounded domain (perforated with a large number of small holes) is modeled through the Navier-Stokes equations with non-homogeneous Dirichlet boundary data satisfying the general outflow condition. Under either a symmetry assumption on the data or under a smallness condition on each of the boundary fluxes (therefore, no constraints on the magnitude of the boundary velocity are imposed), we apply the classical energy method in homogenization theory and study the asymptotic behavior of the solutions to this system as the size of the perforations goes to zero: it is shown that the effective equations remain unmodified in the limit. The main novelty of the present work lies in the obtainment of the required uniform bounds, which are achieved by a contradiction argument based on Bernoulli's law for solutions of the stationary Euler equations.
\end{abstract}	

\section{Introduction and presentation of the problem}
Let $\Omega \subset \mathbb{R}^2$ be an open, bounded and multiply-connected domain of class $\mathcal{C}^{2}$, expressed as
$$
\Omega \doteq \Omega_{0} \setminus \bigcup_{i=1}^{M} \overline{\Omega_{i}} \, ,
$$
for some integer $M \geq 1$ and a collection $\Omega_{0}, \Omega_{1},...,\Omega_{M} \subset \mathbb{R}^{2}$ of open, bounded and simply-connected sets having a $\mathcal{C}^{2}$-boundary such that
$$
\bigcup_{i=1}^{M} \overline{\Omega_{i}} \subset \Omega_{0} \qquad \text{and} \qquad \partial \Omega_{i} \cap \partial \Omega_{j} = \emptyset \quad \forall i,j \in \{ 1,...,M \} \, , \ i \neq j \, .
$$
Then, the boundary $\partial \Omega$ of $\Omega$ consists of $M + 1$ connected components, and is given by
$$
\partial \Omega = \bigcup_{i=0}^{M} \partial\Omega_{i} \, .
$$
For any $\xi \in \mathbb{R}^{2}$ and $r > 0$ we denote by $D(\xi, r) \subset \mathbb{R}^{2}$ the open disk of radius $r$ with center at $\xi$. Let $( K_{n} )_{n \in \mathbb{N}}$ be a sequence of open, bounded and simply connected sets with a $\mathcal{C}^{2}$-boundary such that $(0,0) \in K_{n}$, for every $n \in \mathbb{N}$, and
$$
\sup_{n \in \mathbb{N}} | K_{n} | < \infty \, .
$$
Take $\varepsilon_{*} \in (0,1)$ such that $\varepsilon_{*} | K_{n} | < | \Omega |$, $\forall n \in \mathbb{N}$. Given $\alpha > 2$, we define the quantities
\begin{equation} \label{sigma_eps}
a_\varepsilon \doteq e^{-\frac{1}{\varepsilon^{\alpha}}} \qquad \text{and} \qquad \sigma_\varepsilon \doteq \varepsilon \sqrt{ \left| \log \left( \dfrac{a_\varepsilon}{\varepsilon} \right) \right|} \qquad \forall \varepsilon \in (0,\varepsilon_{*}] \quad \Longrightarrow \quad \lim\limits_{\varepsilon \to 0^{+}} \sigma_{\varepsilon} = +\infty \, .
\end{equation}
Following \cite{nevcasova2023homogenization,nevcasova2022homogenization}, given $\varepsilon \in (0,\varepsilon_{*}]$, suppose that there exist an integer $N(\varepsilon) \geq 1$ and a collection of points $\xi^{\varepsilon}_{1},...,\xi^{\varepsilon}_{N(\varepsilon)} \in \mathbb{R}^{2}$ such that
\begin{equation} \label{perforation}
\begin{aligned}
&  \xi^{\varepsilon}_{n} + a_\varepsilon \overline{K_{n}} \subset D(\xi^{\varepsilon}_{n}, \delta_{0} a_\varepsilon) \subset D(\xi^{\varepsilon}_{n}, \delta_{1} \varepsilon) \subset \Omega \qquad \forall n \in \{1,...,N(\varepsilon)\} \,, \\[6pt]
& \partial D \left(\xi^{\varepsilon}_{n}, \delta_{1} \varepsilon \right) \cap \partial D \left(\xi^{\varepsilon}_{m}, \delta_{1} \varepsilon \right) = \emptyset \qquad \forall n,m \in \{1,...,N(\varepsilon)\} \,, \ n \neq m \, , \\[6pt]
& \partial D \left(\xi^{\varepsilon}_{n}, \delta_{1} \varepsilon \right) \cap \partial \Omega = \emptyset \qquad \forall n \in \{1,...,N(\varepsilon)\} \, ,
\end{aligned}
\end{equation}
for some constants $0 < \delta_{0} < \delta_{1}$ that are independent of $\varepsilon \in (0,\varepsilon_{*}]$. Setting $K^{\varepsilon}_{n} \doteq \xi^{\varepsilon}_{n} + a_\varepsilon K_{n}$ for every $n \in \{1,...,N(\varepsilon)\}$, we will refer to the family $\{ K^{\varepsilon}_{n} \}^{N(\varepsilon)}_{n=1}$ satisfying \eqref{perforation} as the \textit{solid obstacles}, while
\begin{equation} \label{perfordomain}
\Omega_{\varepsilon} \doteq \Omega \setminus \overline{K_{\varepsilon}} \doteq \Omega \setminus \bigcup^{N(\varepsilon)}_{n=1} \overline{K^{\varepsilon}_{n}} \, ,
\end{equation}
represents the \textit{perforated fluid domain} at the $\varepsilon$-level. We emphasize that, given $\varepsilon \in (0,\varepsilon_{*}]$, the family of obstacles $\{ K^{\varepsilon}_{n} \}^{N(\varepsilon)}_{n=1}$ is built in such a way that the \textit{size} of each solid is proportional to $a_\varepsilon$, while the mutual distance between any two consecutive holes is proportional to $\varepsilon$. Moreover, since we only consider those obstacles that are \textit{strictly} contained in $\Omega$ (in the sense of \eqref{perforation}$_3$), the following bound on the number $N(\varepsilon)$ holds:
\begin{equation} \label{number0}
N(\varepsilon) \leq \dfrac{| \Omega |}{\pi \delta_{1}^{2} \varepsilon^{2}} \, .
\end{equation}
\vspace*{-2mm}
\begin{figure}[H]
	\begin{center}
		\includegraphics[scale=0.68]{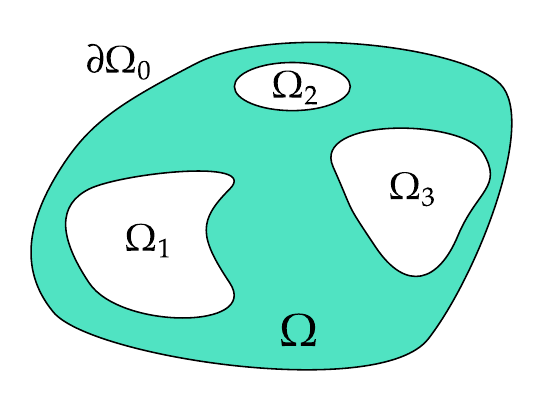}\hspace{2.5cm}
		\includegraphics[scale=0.68]{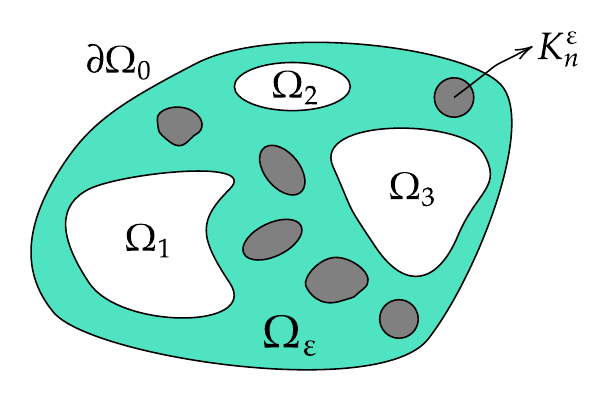}
	\end{center}
	\vspace*{-4mm}
	\caption{Representation of the domain $\Omega$ (left) and the perforated domain $\Omega_{\varepsilon}$ (right).}\label{dom0}
\end{figure}
\noindent
Notice, however, that the solids $\{ K^{\varepsilon}_{n} \}^{N(\varepsilon)}_{n=1}$ may have different shapes and that they are not necessarily periodically distributed in $\Omega$ (compare this with Subsection \ref{perperf}). We decompose the boundary of $\Omega_{\varepsilon}$ as
$$
\partial \Omega_{\varepsilon} = \partial \Omega \cup \partial K_{\varepsilon} = \partial \Omega \cup \bigcup^{N(\varepsilon)}_{n=1} \partial K^{\varepsilon}_{n} \, .
$$
The outward unit normal to $\partial \Omega_{\varepsilon}$ is denoted by $\nu$ (with some abuse of notation, as such vector also depends on $\varepsilon$). Given $\varepsilon \in (0,\varepsilon_{*}]$, we analyze the steady motion of a viscous incompressible fluid (having a constant kinematic viscosity $\eta > 0$) along $\Omega_{\varepsilon}$, which is characterized by its velocity vector field $u_{\varepsilon} : \Omega_{\varepsilon} \longrightarrow \mathbb{R}^2$ and its scalar pressure $p_{\varepsilon} : \Omega_{\varepsilon} \longrightarrow \mathbb{R}$, under the action of an external force $f : \Omega \longrightarrow \mathbb{R}^2$ and a boundary velocity $v_{*} : \partial \Omega \longrightarrow \mathbb{R}^2$ satisfying the compatibility condition
\begin{equation} \label{gof}
\int_{\partial \Omega} v_{*} \cdot \nu = 0 \, .
\end{equation}
Such stationary motion is modeled through the following boundary-value problem (with non-homogeneous Dirichlet boundary conditions) associated to the steady-state Navier-Stokes equations in $\Omega_{\varepsilon}$:
\begin{equation}\label{nsstokes0}
\left\{
\begin{aligned}
& -\eta\Delta u_{\varepsilon}+(u_{\varepsilon}\cdot\nabla)u_{\varepsilon}+\nabla p_{\varepsilon}=f \, , \quad  \nabla\cdot u_{\varepsilon}=0 \ \ \mbox{ in } \ \ \Omega_{\varepsilon} \, , \\[4pt]
& u_{\varepsilon}=v_{*} \ \ \mbox{ on } \ \ \partial \Omega \, , \quad u_{\varepsilon}=0 \ \ \mbox{ on } \ \ \partial K_{\varepsilon} \, .
\end{aligned}
\right.
\end{equation}

\subsection{Symmetric domain with periodic perforations} \label{perperf}

Special attention will be devoted to the case when each of the domains $\Omega_{0}, \Omega_{1},...,\Omega_{M}$ intersects the $x$-axis and is \textit{symmetric with respect to the $x$-axis}, meaning that
$$
(x,y) \in \Omega_{j} \quad \Longleftrightarrow \quad (x,-y) \in \Omega_{j} \qquad \forall j \in \{0,...,M \} \, .
$$
Let $K \subset \mathbb{R}^{2}$ be an open, bounded and simply connected set with a $\mathcal{C}^{2}$-boundary such that $(0,0) \in K$. Suppose there exists $\lambda_{0} \in (0,1)$ for which
$$
\overline{K} \subset D((0,0), \lambda_{0}) \subset \overline{D((0,0), \lambda_{0})} \subset (0,1) \times (0,1) \, ,
$$
and take $\varepsilon_{*} \in (0,1)$ verifying $\varepsilon_{*} | K | < | \Omega |$. Given $\varepsilon \in (0,\varepsilon_{*}]$, for every $m \in \mathbb{Z}^{2}$ we then have
\begin{equation} \label{perfordomainsim0}
\overline{K^{\varepsilon}_{m}} \doteq \varepsilon m + a_\varepsilon \overline{K} \subset D(\varepsilon m, \lambda_{0} a_\varepsilon) \subset \overline{D(\varepsilon m, \lambda_{0} a_\varepsilon)} \subset \varepsilon m + (0, \varepsilon) \times (0, \varepsilon) \, .
\end{equation}
Following \cite{allaire1, conca1985application, feireisl2015homogenization}, the \textit{periodically perforated fluid domain} at the $\varepsilon$-level is defined as
\begin{equation} \label{perfordomainsim}
	\Omega_{\varepsilon} \doteq \Omega \setminus \overline{K_{\varepsilon}} \doteq \Omega \setminus \bigcup_{m \in \mathcal{M}_{\varepsilon}} \overline{K^{\varepsilon}_{m}} \qquad \text{where} \qquad \mathcal{M}_{\varepsilon} \doteq \{ m \in \mathbb{Z}^{2} \ | \ \varepsilon m + (0, \varepsilon) \times (0, \varepsilon) \subset \Omega \} \, .
\end{equation}
Therefore, $\Omega_{\varepsilon}$ is obtained by removing from $\Omega$ all translations (through vectors $\varepsilon m$, with $m \in \mathbb{Z}^{2}$) of the representative hole $a_\varepsilon \overline{K}$ that are \textit{strictly} contained in $\Omega$ (see the last inclusion in \eqref{perfordomainsim0}). Thus,
$$
| \mathcal{M}_{\varepsilon} | \leq \dfrac{| \Omega |}{\varepsilon^{2}} \, .
$$
\vspace*{-5mm}
\begin{figure}[H]
	\begin{center}
			\includegraphics[scale=0.66]{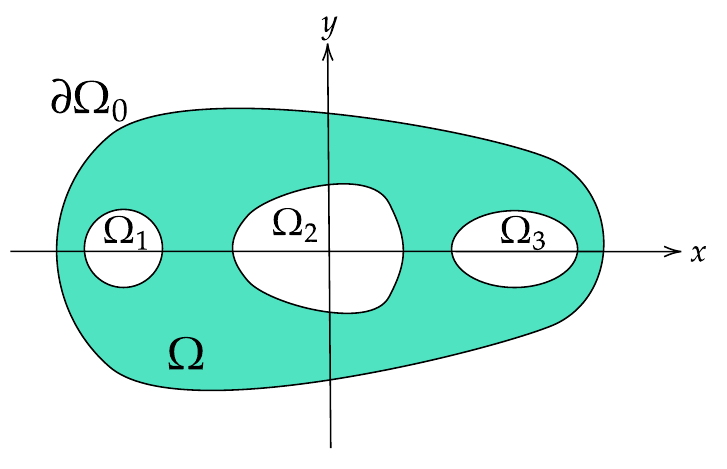}\hspace{1cm}
		\includegraphics[scale=0.66]{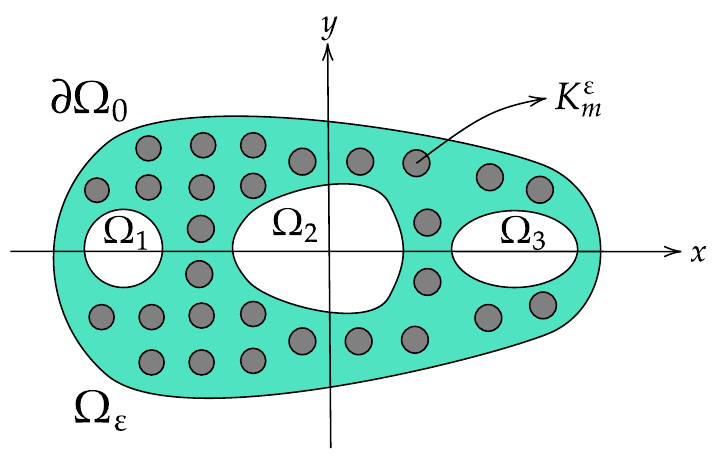}
	\end{center}
	\vspace*{-5mm}
	\caption{The symmetric domain $\Omega$ (left) and the symmetric perforated domain $\Omega_{\varepsilon}$ (right).}\label{dom1}
\end{figure}
\noindent
In particular, $\Omega_{\varepsilon}$ is also symmetric with respect to the $x$-axis, for every $\varepsilon \in (0,\varepsilon_{*}]$.
\par
In order to determine the existence of symmetric (with respect the $x$-axis) generalized solutions to problem \eqref{nsstokes0}, the following definition is given (see \cite[Introduction]{amick1984existence} and \cite{morimoto2007remark,sazonov1993existence}):

\begin{definition} \hlabel{def:axisim}
Let $Q \subset \mathbb{R}^{2}$ be a symmetric domain with respect to the $x$-axis.
	\begin{itemize}[leftmargin=*]
		\item[-] We say that a scalar function $g : Q \longrightarrow \mathbb{R}$ is \textbf{symmetric} (with respect to the $x$-axis) if
		$$
		g(x,-y) = g(x,y) \qquad \forall (x,y) \in Q \, .
		$$
        \item[-] We say that a vector field $G = (G_{1},G_{2}) : Q \longrightarrow \mathbb{R}^{2}$ is \textbf{symmetric} (with respect to the $x$-axis) if
		$$
		G_{1}(x,-y) = G_{1}(x,y) \qquad \text{and} \qquad G_{2}(x,-y) = -G_{2}(x,y) \qquad \forall (x,y) \in Q \, .
		$$
	\end{itemize}
\end{definition}

System \eqref{nsstokes0} constitutes a particular instance of the celebrated \textit{Leray flux problem} \cite{korobkov2014flux, korobkov2016leray} in the context of homogenization theory. For a fixed $\varepsilon \in (0,\varepsilon_{*}]$, the weak solvability of \eqref{nsstokes0} under the \textbf{general outflow condition} \eqref{gof} (which, physically, allows for the presence of interior sources and sinks in $\Omega_{\varepsilon}$) remained one of the most outstanding open problems in Mathematical Fluid Mechanics between the years 1933 and 2015, listed among the \textit{Eleven Great Problems of Mathematical Hydrodynamics} by Yudovich \cite{yudovich2003eleven}. In fact, assuming the \textbf{stringent outflow condition}
$$
\int_{\partial \Omega_{j}} v_{*} \cdot \nu = 0 \qquad \forall j \in \{0,...,M \} \, ,
$$
which rules out the presence of sources and sinks inside $\Omega_{\varepsilon}$, Leray \cite{leray1933etude} showed in 1933 the existence of a generalized solution to \eqref{nsstokes0} employing a fixed-point technique, nowadays known as the \textit{Leray-Schauder Principle} \cite[Chapter 6]{zeidler2013nonlinear}, combined with a \textit{reductio ad absurdum} argument. Independently from each other, Amick \cite{amick1984existence} and Sazonov \cite{sazonov1993existence} recovered the result of Leray under the general outflow condition \eqref{gof}, but imposing a \textit{symmetry assumption} on both $\Omega_{\varepsilon}$ and the data of the problem (external force and boundary velocity): while Amick extended the original method of Leray, Sazonov built a symmetric \textit{solenoidal extension} of the boundary velocity satisfying the \textit{Leray-Hopf inequality} \cite{galdi1991existence}. Relaxing this symmetry restriction, Neustupa \cite{neustupa2010new} reached the same outcome under a smallness condition on each of the boundary fluxes. The fully two-dimensional general result was proved in 2015 by Korobkov, Pileckas \& Russo \cite{korobkov2015solution}, where the uniform bounds required by the Leray-Schauder Principle were achieved by a contradiction argument based on Bernoulli's law \cite{kapitanskii1983spaces} for solutions of the stationary Euler equations and a generalization of the Morse-Sard Theorem for Sobolev functions given by Bourgain, Kristensen \& Korobkov in \cite{bourgain2013morse}. Nevertheless, as far as our knowledge goes, the corresponding homogenization limit associated to Leray's flux problem has not been tackled before, and constitutes the core of the present article. Precisely, our goal is to study the asymptotic behavior of the solutions of problem \eqref{nsstokes0} as $\varepsilon \to 0^{+}$ in two different settings:
\begin{itemize}[leftmargin=*]
\item[-] in the general perforated domain \eqref{perfordomain} under a smallness condition on each of the boundary fluxes;
\item[-] in the symmetric perforated domain \eqref{perfordomainsim} under a symmetry condition on the data of the problem.
\end{itemize}
Notice, therefore, that neither of the above settings imposes a size constraint on the boundary velocity, see Remark \ref{rem1}. To achieve our goal, in Section \ref{epslevelsec} we derive uniform $\varepsilon$-independent bounds for the solutions of \eqref{nsstokes0}. The proofs of Theorems \ref{epslevel}-\ref{epslevelsim}, the most involved in this work, adapt the contradiction argument of \cite{korobkov2015solution} previously described and, additionally, employ several properties of the \textit{relative capacity} of the perforations inside $\Omega$ (see Lemma \ref{refcapacitypro}) and a uniform solenoidal extension of the boundary velocity, see Lemma \ref{bogtype} and its symmetric counterpart Corollary \ref{bogtypesim}; in turn, Lemma \ref{bogtype} relies on the $\varepsilon$-uniform boundedness of the Bogovskii operator in these perforated domains, as shown by Ne{\v{c}}asov{\'a} et al. \cite{nevcasova2023homogenization,nevcasova2022homogenization} (see \cite{diening2017inverse} for the corresponding result in the three-dimensional case). Subsequently, applying the classical energy method in homogenization theory \cite[Appendix]{sanchez1980non}, in Section \ref{energymethod} we show that, as $\varepsilon \to 0^{+}$, the effective or homogenized equation remains unchanged in the limit: up to the extraction of a subsequence, the sequence of solutions (indexed by the parameter $\varepsilon$) of \eqref{nsstokes0} converges \textit{strongly} (in a sense made precise in Theorems \ref{effectiveq1}-\ref{effectiveq2}) to a weak solution of problem \eqref{nsstokes0} in $\Omega$ as $\varepsilon \to 0^{+}$. We point out that the statements of Theorems \ref{effectiveq1}-\ref{effectiveq2} reflect other existing results in the homogenization literature in the regime of \textit{very} tiny holes \cite{allaire2,feireisl2015homogenization, lu2020homogenization,lu2021homogenization,lu2018homogenization, lu2023homogenization}. Finally, Section \ref{remarksfulcase} provides some remarks concerning the possibility of recovering the results contained in this article in the fully two-dimensional case, that is, dropping both the symmetry condition and the smallness assumption on the boundary fluxes.

\newpage
\section{Boundary-value problem at the $\varepsilon$-level: uniform bounds} \label{epslevelsec}
Let $\varepsilon \in I_{*}$ be a fixed parameter, with $I_{*} \doteq (0,\varepsilon_{*}]$. Many of the results contained in the present article exploit the concept of \textit{relative capacity} of $K_{\varepsilon}$ with respect to $\Omega$, defined as
\begin{equation} \label{relcap}
	\mbox{Cap}_{\Omega}(K_{\varepsilon}) \doteq \min_{v\in H^{1}_{0}(\Omega)} \left\{  \int_{\Omega} |\nabla v|^2  \ \Big| \ v=1 \  \text{ in } \ \overline{K_{\varepsilon}} \right\} \, .
\end{equation}
The \textit{relative capacity potential} of $K_{\varepsilon}$ with respect to $\Omega$, that is, the scalar function $\phi_{\varepsilon} \in H_{0}^{1}(\Omega)$ achieving the minimum in \eqref{relcap}, satisfies
\begin{equation} \label{cap1}
	\Delta \phi_{\varepsilon}=0 \ \text{ in } \ \Omega_{\varepsilon} \, , \qquad \phi_{\varepsilon} = 0 \ \text{ on } \ \partial \Omega\, , \qquad
	\phi_{\varepsilon} = 1 \ \text{ in } \ \overline{K_{\varepsilon}} \, , \qquad\mbox{Cap}_{\Omega}(K_{\varepsilon})=\|\nabla \phi_{\varepsilon} \|^2_{L^{2}(\Omega)} \, ,
\end{equation}
see \cite[Chapter 2]{maz2013sobolev} for more details. Further essential properties of the relative capacity potential are collected in the following result, in the spirit of  \cite[Proposition 4.3]{allaire3} and the examples of \cite[Section 2]{cioranescu2018strange}:
\begin{lemma} \label{refcapacitypro}
	Let $\Omega_{\varepsilon}$ be as in \eqref{perfordomain} and $\phi_{\varepsilon} \in H_{0}^{1}(\Omega)$ be the function satisfying \eqref{cap1}. Then, $\phi_{\varepsilon} \in H^{2}(\Omega_{\varepsilon})$ and the following estimates hold
	\begin{equation} \label{cap2}
		\|1-\phi_{\varepsilon} \|_{L^{\infty}(\Omega_{\varepsilon})} \leq C_{*} \, \qquad \text{and} \qquad \| \phi_{\varepsilon} \|_{L^{2}(\Omega_{\varepsilon})} + \|\nabla \phi_{\varepsilon} \|_{L^{2}(\Omega_{\varepsilon})} \leq \dfrac{C_{*}}{\sigma_{\varepsilon}} \, ,
	\end{equation}
	for some constant $C_{*} > 0$ that depends on $\Omega$ and $\{ \delta_{0}, \delta_{1} \}$, but is independent of $\varepsilon \in I_{*}$.
\end{lemma}
\noindent
\begin{proof}
	In what follows, $C > 0$ will always denote a generic constant that depends on $\Omega$ and $\{ \delta_{0}, \delta_{1} \}$ (independently of $\varepsilon \in I_{*}$), but that may change from line to line.
	\newline	
	Since $\Omega_{\varepsilon}$ has a boundary of class $\mathcal{C}^{2}$, standard elliptic regularity arguments show that $\phi_{\varepsilon} \in H^{2}(\Omega_{\varepsilon})$. The first estimate in \eqref{cap2} follows directly from the Maximum Principle. Concerning the second estimate in \eqref{cap2}, given $n \in \{1,...,N(\varepsilon)\}$, consider the function $\varphi^{\varepsilon}_{n} : \Omega \longrightarrow \mathbb{R}$ defined by
	$$
	\varphi^{\varepsilon}_{n}(\xi) =
	\begin{cases}
		1 & \quad \text{if} \ \ 0 \leq | \xi - \xi^{\varepsilon}_{n} | \leq \delta_{0} a_\varepsilon  \, ,\\[3pt]
		\dfrac{\log(\delta_{1} \varepsilon) - \log(| \xi - \xi^{\varepsilon}_{n} |)}{\log(\delta_{1} \varepsilon) - \log(\delta_{0} a_\varepsilon)} & \quad  \text{if} \ \ \delta_{0} a_\varepsilon < | \xi - \xi^{\varepsilon}_{n} | \leq \delta_{1} \varepsilon \, ,\\[3pt]
		0 & \quad  \text{if} \ \ | \xi - \xi^{\varepsilon}_{n} | > \delta_{1} \varepsilon \, ,
	\end{cases}
	$$
	so that, by the assumptions in \eqref{perforation}$_1$, $\varphi^{\varepsilon}_{n} \in H^{1}_{0}(\Omega)$ and $\varphi^{\varepsilon}_{n} = 1$ in $\overline{D(\xi^{\varepsilon}_{n}, \delta_{0} a_\varepsilon)}$. As the relative capacity is an outer measure and is increasing with respect to domain inclusion (see \cite[Section 2.2]{maz2013sobolev}), we get
	$$
	\mbox{Cap}_{\Omega}(K_{\varepsilon}) \leq \sum_{n=1}^{N(\varepsilon)} \mbox{Cap}_{\Omega}(K^{\varepsilon}_{n}) \leq \sum_{n=1}^{N(\varepsilon)} \mbox{Cap}_{\Omega}(D(\xi^{\varepsilon}_{n}, \delta_{0} a_\varepsilon)) \leq \sum_{n=1}^{N(\varepsilon)} \int_{\Omega} |\nabla \varphi^{\varepsilon}_{n}|^2 \leq \dfrac{C}{\sigma_{\varepsilon}^{2}} \, ,
	$$
	so that
	$$
	\| \nabla \phi_{\varepsilon} \|_{L^{2}(\Omega)} \leq \dfrac{C}{\sigma_{\varepsilon}} \, .
	$$
	Since $\phi_{\varepsilon} \in H_{0}^{1}(\Omega)$, an application of Poincaré's inequality in $\Omega$ allows us to conclude the proof.
\end{proof}

Let $Q \subset \mathbb{R}^{2}$ be any bounded Lipschitz domain, and consider the space of square-integrable functions in $Q$ having zero mean value:
$$
	L^2_0(Q) \doteq \left\{ g\in L^2(Q) \ \Big | \ \int_Q g = 0 \right\} \, .
$$
Another essential preliminary result concerns the construction of a uniform \textit{solenoidal extension} of the boundary velocity. To fix notation, given any function (scalar or vector) $\psi \in L^{1}(\Omega_{\varepsilon})$, hereafter we will denote by $\widetilde{\psi} \in L^{1}(\Omega)$ the function defined by
$$
\widetilde{\psi} \doteq
\begin{cases}
	\psi & \quad \text{in} \ \ \Omega_{\varepsilon} \, ,\\[3pt]
	0 & \quad  \text{in} \ \ \overline{K_{\varepsilon}}  \, .
\end{cases}
$$
We then prove:
\begin{lemma} \label{bogtype}
	Let $\Omega_{\varepsilon}$ be as in \eqref{perfordomain} and $v_{*} \in H^{3/2}(\partial \Omega)$ satisfying \eqref{gof}. There exists a vector field $J_{\varepsilon} \in H^{1}(\Omega_{\varepsilon})$ such that
	\begin{equation} \label{vecje}
		\left\{
		\begin{aligned}
			& \nabla \cdot J_{\varepsilon} = 0 \ \ \mbox{in} \ \ \Omega_{\varepsilon} \, ; \qquad J_{\varepsilon} = v_{*} \ \ \mbox{on} \ \ \partial \Omega \, ; \\[3pt]
			& J_{\varepsilon} = 0 \ \ \mbox{on} \ \ \partial K_{\varepsilon} \, ; \qquad \| \nabla J_{\varepsilon} \|_{L^{2}(\Omega_{\varepsilon})} \leq C_{*} \| v_{*} \|_{H^{3/2}(\partial \Omega)} \, ,
		\end{aligned}
		\right.
	\end{equation}
	for some constant $C_{*} > 0$ that depends on $\Omega$ and $\{ \delta_{0}, \delta_{1} \}$, but is independent of $\varepsilon \in I_{*}$. In particular, there exists a vector field $\widehat{J} \in H^{1}(\Omega)$ such that
$$
\nabla \cdot \widehat{J} = 0 \ \ \mbox{in} \ \ \Omega \qquad \text{and} \qquad \widehat{J} = v_{*} \ \ \mbox{on} \ \ \partial \Omega \, ,
$$
and a (not relabeled) subsequence $(\widetilde{J_\varepsilon})_{\varepsilon \in I_{*}} \subset H^{1}(\Omega)$ for which
\begin{equation} \label{jepsconv}
\widetilde{J_\varepsilon} \rightharpoonup \widehat{J} \ \ \ \text{weakly in} \ H^{1}(\Omega) \qquad \text{and} \qquad \widetilde{J_\varepsilon} \to \widehat{J} \ \ \ \text{strongly in} \ L^{4}(\Omega) \ \ \ \text{as} \ \ \varepsilon \to 0^{+} \, .
\end{equation}
\end{lemma}
\noindent
\begin{proof}
In what follows, $C > 0$ will always denote a generic constant that depends on $\Omega$ and $\{ \delta_{0}, \delta_{1} \}$ (independently of $\varepsilon \in I_{*}$), but that may change from line to line.
\newline
Since $v_{*} \in H^{3/2}(\partial \Omega)$, there exists a vector field $V_{*} \in H^{2}(\Omega)$ such that
\begin{equation} \label{extens1}
V_{*} = v_{*} \quad \text{on} \ \ \partial \Omega \qquad \text{and} \qquad \| V_{*} \|_{H^{2}(\Omega)} \leq C \| v_{*} \|_{H^{3/2}(\partial \Omega)} \, ,
\end{equation}
see \cite[Teorema 1.I]{gagliardo} or \cite[Theorem II.4.3]{galdi2011introduction} as well. Let $\phi_{\varepsilon} \in H^{2}(\Omega_{\varepsilon}) \cap H_{0}^{1}(\Omega)$ be the relative capacity potential of $K_{\varepsilon}$ with respect to $\Omega$, see Lemma \ref{refcapacitypro}. Then, the vector field $X_{\varepsilon} \doteq (1-\phi_{\varepsilon})V_{*} \in H^{1}(\Omega_{\varepsilon})$ is such that $X_{\varepsilon} = 0$ on $\partial K_{\varepsilon}$ and $X_{\varepsilon} = v_{*}$ on $\partial \Omega$. Moreover
$$
\nabla X_{\varepsilon} = -\nabla \phi_{\varepsilon} \otimes V_{*} + (1 - \phi_{\varepsilon}) \nabla V_{*} \quad \text{in} \ \ \Omega_{\varepsilon} \, ,
$$
and therefore, in view of the Sobolev embedding $H^{2}(\Omega) \subset \mathcal{C}(\overline{\Omega})$ and \eqref{cap2}, we have
\begin{equation} \label{extens2}
\begin{aligned}
\| \nabla X_{\varepsilon} \|_{L^{2}(\Omega_{\varepsilon})} \leq \| \nabla \phi_{\varepsilon} \|_{L^{2}(\Omega_{\varepsilon})} \| V_{*} \|_{L^{\infty}(\Omega)} + \| 1 - \phi_{\varepsilon} \|_{L^{\infty}(\Omega_{\varepsilon})} \| \nabla V_{*} \|_{L^{2}(\Omega)} \leq C \| V_{*} \|_{H^{2}(\Omega)} \, .
\end{aligned}
\end{equation}
Furthermore, from the Divergence Theorem and \eqref{gof} we have
$$
\int_{\Omega_{\varepsilon}} \nabla \cdot X_{\varepsilon} = \int_{\partial \Omega_{\varepsilon}} X_{\varepsilon} \cdot \nu = \int_{\partial \Omega} v_{*} \cdot \nu = 0 \, ,
$$
that is, $\nabla \cdot X_{\varepsilon} \in L^{2}_{0}(\Omega_{\varepsilon})$. Therefore, there exists a vector field $Y_{\varepsilon} \in H^{1}_{0}(\Omega_{\varepsilon})$ such that
\begin{equation} \label{bogn0}
	\nabla \cdot Y_{\varepsilon} = - \nabla \cdot X_{\varepsilon} \quad \text{in} \quad \Omega_{\varepsilon} \qquad \text{and} \qquad  \| \nabla Y_{\varepsilon} \|_{L^{2}(\Omega_{\varepsilon})} \leq C_{B}(\Omega_{\varepsilon}) \| \nabla \cdot X_{\varepsilon} \|_{L^{2}(\Omega_{\varepsilon})}  \, ,
\end{equation}
see \cite{bogovskii1979solution}. From \cite[Proposition 1.7]{nevcasova2022homogenization} we know that $C_{B}(\Omega_{\varepsilon})$ (the so-called \textit{Bogovskii constant} of $\Omega_{\varepsilon}$) admits the uniform bound
\begin{equation} \label{boguni0}
	C_{B}(\Omega_{\varepsilon}) \leq C \, .
\end{equation}
We set $J_{\varepsilon} \doteq X_{\varepsilon} + Y_{\varepsilon}$ which, in view of \eqref{extens1}-\eqref{extens2}-\eqref{bogn0}-\eqref{boguni0}, is an element of $H^{1}(\Omega_{\varepsilon})$ satisfying \eqref{vecje}. Notice that $\widetilde{J_{\varepsilon}} \in H^{1}(\Omega)$ is divergence-free separately in $\Omega_{\varepsilon}$ and $K_{\varepsilon}$. Moreover,
$$
\| \nabla \widetilde{J_{\varepsilon}} \|_{L^{2}(\Omega)} = \| \nabla J_{\varepsilon} \|_{L^{2}(\Omega_{\varepsilon})} \, .
$$
The bound in \eqref{vecje}$_2$ ensures the existence of $\widehat{J} \in H^{1}(\Omega)$ such that the convergences
\begin{equation} \label{convergencesiol}
\widetilde{J_\varepsilon} \rightharpoonup \widehat{J} \ \ \ \text{weakly in} \ H^{1}(\Omega) \, ; \quad \widetilde{J_\varepsilon} \to \widehat{J} \ \ \ \text{strongly in} \ L^{4}(\Omega) \, ; \quad \widetilde{J_\varepsilon} \to \widehat{J} \ \ \ \text{strongly in} \ L^{2}(\partial \Omega)
\end{equation}
hold as $\varepsilon \to 0^{+}$, along a subsequence that is not being relabeled, see also \cite[Theorem 6.2]{necas2011direct}. Given any scalar function $\phi \in \mathcal{C}_{0}^{\infty}(\Omega)$, an integration by parts and the divergence-free condition in \eqref{vecje}$_1$ imply
$$
\int_{\Omega} \widetilde{J_{\varepsilon}} \cdot \nabla \phi = \int_{\Omega_{\varepsilon}} J_{\varepsilon} \cdot \nabla \phi = \int_{\partial \Omega_{\varepsilon}} \phi (J_{\varepsilon} \cdot \nu) = 0 \qquad \forall \varepsilon \in I_{*} \, ,
$$
since $J_{\varepsilon}$ vanishes on $\partial K_{\varepsilon}$ and so does $\phi$ on $\partial \Omega$. Then, the weak convergence in \eqref{convergencesiol} yields
$$
\int_{\Omega} \widehat{J} \cdot \nabla \phi = - \int_{\Omega} \phi (\nabla \cdot \widehat{J} ) = 0 \qquad \forall \phi \in \mathcal{C}_{0}^{\infty}(\Omega; \mathbb{R}) \, ,
$$
thus proving that $\nabla \cdot \widehat{J} = 0$ almost everywhere in $\Omega$. Finally, the strong convergence in $L^{2}(\partial \Omega)$ in \eqref{convergencesiol} also guarantees that $\widehat{J} = v_{*}$ almost everywhere on $\partial \Omega$.
\end{proof}

We introduce the functional spaces (of vector fields) that will be employed hereafter:
$$
H_{\sigma}^{1}(\Omega_{\varepsilon}) \doteq \left\lbrace v \in H^{1}(\Omega_{\varepsilon}) \ | \ \nabla \cdot v=0 \ \ \mbox{in} \ \ \Omega_{\varepsilon} \right\rbrace \qquad \text{and} \qquad H^{1}_{0, \sigma}(\Omega_{\varepsilon}) \doteq \left\lbrace v \in H_{0}^{1}(\Omega_{\varepsilon}) \ | \ \nabla \cdot v=0 \ \ \mbox{in} \ \ \Omega_{\varepsilon} \right\rbrace \, ,
$$
which are Hilbert spaces if endowed with the scalar product of $H^{1}(\Omega_{\varepsilon}; \mathbb{R}^2)$. We recall the standard definition for the weak solutions of problem \eqref{nsstokes0}:
\begin{definition}\label{weaksolution}
	A vector field $u \in H_{\sigma}^{1}(\Omega_{\varepsilon})$ is called a \textbf{weak solution} of \eqref{nsstokes0} if $u - J_{\varepsilon} \in H_{0,\sigma}^{1}(\Omega_{\varepsilon})$ and
	$$
	\eta \int_{\Omega_{\varepsilon}} \nabla u \cdot \nabla \varphi + \int_{\Omega_{\varepsilon}} (u \cdot \nabla)u \cdot \varphi  = \int_{\Omega_{\varepsilon}} f \cdot \varphi \qquad \forall \varphi \in H^{1}_{0,\sigma}(\Omega_{\varepsilon}) \, .
	$$
\end{definition}

Given $v_{*} \in H^{3/2}(\partial \Omega)$ satisfying \eqref{gof}, we will denote by
$$
F_{j} \doteq \int_{\partial \Omega_{j}} v_{*} \cdot \nu \qquad \forall j \in \{0,...,M \}
$$
the flux of $v_{*}$ across each of the connected components of $\partial \Omega$. Given $j \in \{1,...,M \}$, let $\chi_{j} \in H_{0}^{1}(\Omega_{0})$ be the \textit{relative capacity potential} of $\Omega_{j}$ with respect to $\Omega_{0}$, that is, the scalar function satisfying
\begin{equation} \label{harcap1}
	\Delta \chi_{j}=0 \ \text{ in } \ \Omega_{0} \setminus \overline{\Omega_{j}} \, , \qquad \chi_{j} = 0 \ \text{ on } \ \partial \Omega_{0} \,, \qquad \chi_{j} = 1 \ \text{ in } \ \overline{\Omega_{j}} \, , 
\end{equation}
and so we denote by
\begin{equation} \label{harcap2}
\text{Cap}_{\Omega_{0}}(\Omega_{j}) \doteq \int_{\partial \Omega_{j}} \dfrac{\partial \chi_{j}}{\partial \nu}  = \int_{\Omega_{0}} |\nabla \chi_{j}|^2
\end{equation}
the \textit{relative capacity} of $\Omega_{j}$ with respect to $\Omega_{0}$, see \cite[Chapter 2]{maz2013sobolev} for further details. The first main result of this section provides uniform bounds (with respect to $\varepsilon \in I_{*}$) for the solutions of problem \eqref{nsstokes0}.

\begin{theorem} \label{epslevel}
Let $\Omega_{\varepsilon}$ be as in \eqref{perfordomain}. For any $f \in H^{1}(\Omega)$ and $v_{*} \in H^{3/2}(\partial \Omega)$ satisfying \eqref{gof}, there exists at least one weak solution $u_{\varepsilon} \in H^{2}(\Omega_{\varepsilon}) \cap H^{1}_{\sigma}(\Omega_{\varepsilon})$ of problem \eqref{nsstokes0}, and an associated pressure $p_{\varepsilon} \in H^{1}(\Omega_{\varepsilon}) \cap L^{2}_{0}(\Omega_{\varepsilon})$ such that the pair $(u_{\varepsilon},p_{\varepsilon})$ solves \eqref{nsstokes0} in strong form. Moreover, if
\begin{equation} \label{neustupa}
\sum_{j=1}^{M} \dfrac{|F_{j} |}{\textup{Cap}_{\Omega_{0}}(\Omega_{j})} < \eta \, ,
\end{equation}
then the uniform bound
	\begin{equation} \label{uboundpf}
		\sup_{\varepsilon \in I_{*}} \left( \| \nabla u_{\varepsilon} \|_{L^{2}(\Omega_{\varepsilon})} + \| p_{\varepsilon} \|_{L^{2}(\Omega_{\varepsilon})} \right) \leq C_{*} \, ,
	\end{equation}
	holds for some constant $C_{*} > 0$ that depends on $\Omega$, $\eta$, $f$, $v_{*}$ and $\{ \delta_{0}, \delta_{1} \}$, but is independent of $\varepsilon \in I_{*}$.
\end{theorem}
\noindent
\begin{proof}
In what follows, $C > 0$ will always denote a generic constant that depends on $\Omega$, $\eta$ and $\{ \delta_{0}, \delta_{1} \}$ (independently of $\varepsilon \in I_{*}$), but that may change from line to line.
\par
Then, given any $f \in H^{1}(\Omega)$, $v_{*} \in H^{3/2}(\partial \Omega)$ satisfying \eqref{gof} and $\varepsilon \in I_{*}$, \cite[Theorem 1.1]{korobkov2015solution} (see also \cite[Remark 1.1]{korobkov2015solution}) ensures the existence of at least one weak solution $u_{\varepsilon} \in H^{2}(\Omega_{\varepsilon}) \cap H^{1}_{\sigma}(\Omega_{\varepsilon})$ of problem \eqref{nsstokes0} and an associated pressure $p_{\varepsilon} \in H^{1}(\Omega_{\varepsilon}) \cap L^{2}_{0}(\Omega_{\varepsilon})$ such that the pair $(u_{\varepsilon},p_{\varepsilon})$ satisfies \eqref{nsstokes0} in strong form. Let $J_{\varepsilon} \in H^{1}(\Omega_{\varepsilon})$ be the vector field arising from Lemma \eqref{bogtype}, which satisfies \eqref{vecje}. We set $w_{\varepsilon} \doteq u_{\varepsilon} - J_{\varepsilon} \in H^{1}_{0,\sigma}(\Omega_{\varepsilon})$, so that, in view of the bound in \eqref{vecje}$_2$, it suffices to show that
$$
	\sup_{\varepsilon \in I_{*}} \left( \| \nabla w_{\varepsilon} \|_{L^{2}(\Omega_{\varepsilon})} + \| p_{\varepsilon} \|_{L^{2}(\Omega_{\varepsilon})} \right) \leq C \, .
$$
The weak formulation in Definition \eqref{weaksolution} may be then re-written as
\begin{equation} \label{weaksolutionsh}
\begin{aligned}
& \eta \int_{\Omega_{\varepsilon}} \nabla w_{\varepsilon} \cdot \nabla \varphi + \int_{\Omega_{\varepsilon}} (w_{\varepsilon} \cdot \nabla)w_{\varepsilon} \cdot \varphi + \int_{\Omega_{\varepsilon}} (w_{\varepsilon} \cdot \nabla) J_{\varepsilon} \cdot \varphi + \int_{\Omega_{\varepsilon}} (J_{\varepsilon} \cdot \nabla)w_{\varepsilon} \cdot \varphi \\[6pt]
& \hspace{-4mm} = \int_{\Omega_{\varepsilon}} f \cdot \varphi - \eta \int_{\Omega_{\varepsilon}} \nabla J_{\varepsilon} \cdot \nabla \varphi - \int_{\Omega_{\varepsilon}} (J_{\varepsilon} \cdot \nabla)J_{\varepsilon} \cdot \varphi \qquad \forall \varphi \in H^{1}_{0,\sigma}(\Omega_{\varepsilon}) \, .
\end{aligned}
\end{equation}
For later use, we also present the weak formulation of Definition \eqref{weaksolution} incorporating the pressure term:
\begin{equation} \label{weaksolutionshp}
	\begin{aligned}
		& \eta \int_{\Omega_{\varepsilon}} \nabla w_{\varepsilon} \cdot \nabla \varphi + \int_{\Omega_{\varepsilon}} (w_{\varepsilon} \cdot \nabla)w_{\varepsilon} \cdot \varphi + \int_{\Omega_{\varepsilon}} (w_{\varepsilon} \cdot \nabla) J_{\varepsilon} \cdot \varphi + \int_{\Omega_{\varepsilon}} (J_{\varepsilon} \cdot \nabla)w_{\varepsilon} \cdot \varphi - \int_{\Omega_{\varepsilon}} p_{\varepsilon} (\nabla \cdot \varphi) \\[6pt]
		& \hspace{-4mm} = \int_{\Omega_{\varepsilon}} f \cdot \varphi - \eta \int_{\Omega_{\varepsilon}} \nabla J_{\varepsilon} \cdot \nabla \varphi - \int_{\Omega_{\varepsilon}} (J_{\varepsilon} \cdot \nabla)J_{\varepsilon} \cdot \varphi \qquad \forall \varphi \in H^{1}_{0}(\Omega_{\varepsilon}) \, .
	\end{aligned}
\end{equation}
Notice that $\widetilde{p_{\varepsilon}} \in L_{0}^{2}(\Omega)$ and that $\widetilde{w_{\varepsilon}} \in H_{0}^{1}(\Omega)$ is divergence-free separately in $\Omega_{\varepsilon}$ and $K_{\varepsilon}$. Moreover,
\begin{equation}\label{extvel}
	\| \nabla \widetilde{w_{\varepsilon}} \|_{L^{2}(\Omega)} = \| \nabla w_{\varepsilon} \|_{L^{2}(\Omega_{\varepsilon})} \qquad \text{and} \qquad \| \widetilde{p_{\varepsilon}} \|_{L^{2}(\Omega)} = \| p_{\varepsilon} \|_{L^{2}(\Omega_{\varepsilon})} \, .
\end{equation}
We take $\varphi = w_{\varepsilon}$ in \eqref{weaksolutionsh} and integrate by parts, thereby obtaining:
\begin{equation}\label{test2}
\eta \| \nabla w_{\varepsilon} \|^{2}_{L^{2}(\Omega_{\varepsilon})} = - \eta \int_{\Omega} \nabla \widetilde{w_{\varepsilon}} \cdot \nabla \widetilde{J_{\varepsilon}} + \int_{\Omega} (\widetilde{w_{\varepsilon}} \cdot \nabla)\widetilde{w_{\varepsilon}} \cdot \widetilde{J_{\varepsilon}} - \int_{\Omega} (\widetilde{J_{\varepsilon}} \cdot \nabla)\widetilde{J_{\varepsilon}} \cdot \widetilde{w_{\varepsilon}}  + \int_{\Omega} f \cdot \widetilde{w_{\varepsilon}}  \, .
\end{equation}
Now, since $p_{\varepsilon} \in L^{2}_{0}(\Omega_{\varepsilon})$, let $P_{\varepsilon} \in H_{0}^{1}(\Omega_{\varepsilon})$ be a vector field such that
\begin{equation} \label{vecjepre}
	\nabla \cdot P_{\varepsilon} = p_{\varepsilon} \quad \text{in} \quad \Omega_{\varepsilon} \qquad \text{and} \qquad  \| \nabla P_{\varepsilon} \|_{L^{2}(\Omega_{\varepsilon})} \leq C \| p_{\varepsilon} \|_{L^{2}(\Omega_{\varepsilon})} \, ,
\end{equation}
see \eqref{bogn0}-\eqref{boguni0}. We multiply the first identity in \eqref{nsstokes0}$_1$ by $P_{\varepsilon}$ and integrate by parts in $\Omega_{\varepsilon}$ to obtain
$$
\eta \int_{\Omega_{\varepsilon}} \nabla u_{\varepsilon} \cdot \nabla P_{\varepsilon} + \int_{\Omega_{\varepsilon}} (u_{\varepsilon} \cdot \nabla)u_{\varepsilon} \cdot P_{\varepsilon} - \| p_{\varepsilon} \|^{2}_{L^{2}(\Omega_{\varepsilon})} = \int_{\Omega_{\varepsilon}} f \cdot P_{\varepsilon} \, .
$$
Observing that $\widetilde{P_{\varepsilon}} \in H_{0}^{1}(\Omega)$ and that $\| \nabla \widetilde{P_{\varepsilon}} \|_{L^{2}(\Omega)}=\| \nabla P_{\varepsilon} \|_{L^{2}(\Omega_{\varepsilon})}$, we apply H\"older's inequality, the Sobolev and Poincaré inequalities in $\Omega$, and ultimately \eqref{extvel}-\eqref{vecjepre}, in order to estimate
$$
\begin{aligned}
\| p_{\varepsilon} \|^{2}_{L^{2}(\Omega_{\varepsilon})} & = \eta \int_{\Omega_{\varepsilon}} \nabla u_{\varepsilon} \cdot \nabla P_{\varepsilon} + \int_{\Omega} (\widetilde{u_{\varepsilon}} \cdot \nabla)\widetilde{u_{\varepsilon}} \cdot \widetilde{P_{\varepsilon}} - \int_{\Omega} f \cdot \widetilde{P_{\varepsilon}} \\[6pt]
& \leq \eta \| \nabla u_{\varepsilon} \|_{L^{2}(\Omega_{\varepsilon})} \| \nabla P_{\varepsilon} \|_{L^{2}(\Omega_{\varepsilon})} + \| \nabla \widetilde{u_{\varepsilon}} \|_{L^{2}(\Omega)} \| \widetilde{u_{\varepsilon}} \|_{L^{4}(\Omega)} \| \widetilde{P_{\varepsilon}} \|_{L^{4}(\Omega)} + \| f \|_{L^{2}(\Omega)} \| \widetilde{P_{\varepsilon}} \|_{L^{2}(\Omega)} \\[6pt]
& \leq \eta \| \nabla u_{\varepsilon} \|_{L^{2}(\Omega_{\varepsilon})} \| \nabla P_{\varepsilon} \|_{L^{2}(\Omega_{\varepsilon})} + C \| \nabla \widetilde{u_{\varepsilon}} \|^{2}_{L^{2}(\Omega)} \| \nabla \widetilde{P_{\varepsilon}} \|_{L^{2}(\Omega)} + C \| f \|_{L^{2}(\Omega)} \| \nabla \widetilde{P_{\varepsilon}} \|_{L^{2}(\Omega)} \\[6pt]
& \leq C \left( \| \nabla u_{\varepsilon} \|_{L^{2}(\Omega_{\varepsilon})} + \| \nabla u_{\varepsilon} \|^{2}_{L^{2}(\Omega_{\varepsilon})} + \| f \|_{L^{2}(\Omega)} \right) \| \nabla P_{\varepsilon} \|_{L^{2}(\Omega_{\varepsilon})} \\[6pt]
& \leq C \left( 1 + \| \nabla u_{\varepsilon} \|^{2}_{L^{2}(\Omega_{\varepsilon})} + \| f \|_{L^{2}(\Omega)} \right) \| p_{\varepsilon} \|_{L^{2}(\Omega_{\varepsilon})} \, ,
\end{aligned}
$$
so that, after further applying \eqref{vecje}$_2$, we get
\begin{equation} \label{estpresind}
\begin{aligned}
\| p_{\varepsilon} \|_{L^{2}(\Omega_{\varepsilon})} & \leq C \left( 1 + \| \nabla u_{\varepsilon} \|^{2}_{L^{2}(\Omega_{\varepsilon})} + \| f \|_{L^{2}(\Omega)} \right) \\[6pt]
& \leq C \left( 1 + \| \nabla w_{\varepsilon} \|^{2}_{L^{2}(\Omega_{\varepsilon})} + \| v_{*} \|^{2}_{H^{3/2}(\partial \Omega)} +  \| f \|_{L^{2}(\Omega)} \right) \qquad \forall \varepsilon \in I_{*} \, .
\end{aligned}
\end{equation}
By contradiction, suppose now that the norms $\| \nabla w_{\varepsilon} \|_{L^{2}(\Omega_{\varepsilon})}$ are not uniformly bounded with respect to $\varepsilon \in I_{*}$. Then, there must exist a subsequence (not being relabeled) such that
\begin{equation} \label{divergent}
\lim_{\varepsilon \to 0^{+}} \mathcal{Z}_{\varepsilon} = + \infty  \quad \text{with} \ \ \mathcal{Z}_{\varepsilon} \doteq \| \nabla w_{\varepsilon} \|_{L^{2}(\Omega_{\varepsilon})} \quad \forall \varepsilon \in I_{*} \, .
\end{equation}
The estimate \eqref{estpresind} enables us to establish that, along this divergent subsequence \eqref{divergent}, the following sequences are uniformly bounded with respect to $\varepsilon \in I_{*}$:
$$
( \widehat{w_\varepsilon} )_{\varepsilon \in I_{*}} \doteq \left( \dfrac{\widetilde{w_{\varepsilon}}}{\mathcal{Z}_{\varepsilon}} \right)_{\varepsilon \in I_{*}} \subset H_{0}^{1}(\Omega) \qquad \text{and} \qquad ( \widehat{p_\varepsilon} )_{\varepsilon \in I_{*}} \doteq  \left( \dfrac{\widetilde{p_{\varepsilon}}}{\mathcal{Z}^{2}_{\varepsilon}} \right)_{\varepsilon \in I_{*}} \subset L^{2}(\Omega) \, .
$$
Then, there must exist $\widehat{w} \in H_{0}^{1}(\Omega)$ and $\widehat{p} \in L^{2}(\Omega)$ such that the following convergences hold:
\begin{equation} \label{convergencesn11}
\begin{aligned}
\widehat{w_\varepsilon} \rightharpoonup \widehat{w} \ \ \ \text{weakly in} \ H^{1}(\Omega) \, ; \qquad &\widehat{w_\varepsilon} \to \widehat{w} \ \ \ \text{strongly in} \ L^{4}(\Omega) \, ; \qquad \widehat{p_\varepsilon} \rightharpoonup \widehat{p} \ \ \ \text{weakly in} \ L^{2}(\Omega) \, ,
\end{aligned}
\end{equation}
as $\varepsilon \to 0^{+}$, along subsequences that are not being relabeled. Clearly we have $\| \nabla \widehat{w} \|_{L^{2}(\Omega)} \leq 1$, and exactly as in the final part of the proof of Lemma \eqref{bogtype} we can easily deduce that $\nabla \cdot \widehat{w} = 0$ almost everywhere in $\Omega$, that is, $\widehat{w} \in H_{0,\sigma}^{1}(\Omega)$. If we then divide identity \eqref{test2} by $\mathcal{Z}^{2}_{\varepsilon}$, we get
\begin{equation}\label{test3}
	\eta  = - \dfrac{\eta}{\mathcal{Z}_{\varepsilon}} \int_{\Omega} \nabla \widehat{w_{\varepsilon}} \cdot \nabla \widetilde{J_{\varepsilon}} + \int_{\Omega} (\widehat{w_{\varepsilon}} \cdot \nabla)\widehat{w_{\varepsilon}} \cdot \widetilde{J_{\varepsilon}} - \dfrac{1}{\mathcal{Z}_{\varepsilon}} \int_{\Omega} (\widetilde{J_{\varepsilon}} \cdot \nabla)\widetilde{J_{\varepsilon}} \cdot \widehat{w_{\varepsilon}}  + \dfrac{1}{\mathcal{Z}_{\varepsilon}} \int_{\Omega} f \cdot \widehat{w_{\varepsilon}} \qquad \forall \varepsilon \in I_{*} \, ,
\end{equation}
along the subsequences \eqref{convergencesn11}. In order to handle the second term appearing at the right-hand side of \eqref{test3} we firstly notice that
\begin{equation}\label{test4}
\int_{\Omega} (\widehat{w_{\varepsilon}} \cdot \nabla)\widehat{w_{\varepsilon}} \cdot \widetilde{J_{\varepsilon}} = \int_{\Omega} (\widehat{w} \cdot \nabla)\widehat{w_{\varepsilon}} \cdot \widehat{J} + \int_{\Omega} ((\widehat{w_{\varepsilon}}-\widehat{w}) \cdot \nabla)\widehat{w_{\varepsilon}} \cdot \widehat{J} + \int_{\Omega} (\widehat{w_{\varepsilon}} \cdot \nabla)\widehat{w_{\varepsilon}} \cdot (\widetilde{J_{\varepsilon}} - \widehat{J}) \qquad \forall \varepsilon \in I_{*} \, .
\end{equation}
On one hand, the weak convergence in \eqref{convergencesn11} implies
\begin{equation} \label{nsstokeslambda1}
	\lim_{\varepsilon \to 0^{+}} \int_{\Omega} (\widehat{w} \cdot \nabla)\widehat{w_{\varepsilon}} \cdot \widehat{J} = \int_{\Omega} (\widehat{w} \cdot \nabla)\widehat{w} \cdot \widehat{J} \, .
\end{equation}
On the other hand, the H\"older and Sobolev inequalities in $\Omega$, together with the strong convergences in \eqref{jepsconv}-\eqref{convergencesn11} entail that
\begin{equation} \label{bypartsje55}
	\begin{aligned}
		& \left| \int_{\Omega} ((\widehat{w_{\varepsilon}}-\widehat{w}) \cdot \nabla)\widehat{w_{\varepsilon}} \cdot \widehat{J} + \int_{\Omega} (\widehat{w_{\varepsilon}} \cdot \nabla)\widehat{w_{\varepsilon}} \cdot (\widetilde{J_{\varepsilon}} - \widehat{J})  \right| \\[6pt]
		& \hspace{-4.5mm} \leq \left| \int_{\Omega} ((\widehat{w_{\varepsilon}}-\widehat{w}) \cdot \nabla)\widehat{w_{\varepsilon}} \cdot \widehat{J} \right| + \left| \int_{\Omega} (\widehat{w_{\varepsilon}} \cdot \nabla)\widehat{w_{\varepsilon}} \cdot (\widetilde{J_{\varepsilon}} - \widehat{J})  \right| \\[6pt]
		& \hspace{-4.5mm} \leq C \left( \| \widehat{w_{\varepsilon}} - \widehat{w} \|_{L^{4}(\Omega)} + \| \widetilde{J_{\varepsilon}} - \widehat{J} \|_{L^{4}(\Omega)} \right) \to 0 \quad \text{as} \ \ \varepsilon \to 0^{+} \, .
	\end{aligned}
\end{equation}
Inserting \eqref{test4}-\eqref{nsstokeslambda1}-\eqref{bypartsje55} into \eqref{test3} and letting $\varepsilon \to 0^{+}$, again from \eqref{jepsconv}-\eqref{convergencesn11} we deduce
\begin{equation}\label{contradick}
	\eta  =  \int_{\Omega} (\widehat{w} \cdot \nabla)\widehat{w} \cdot \widehat{J} \, .
\end{equation}
A contradiction will be reached in \eqref{contradick} once we prove that
$$
\left| \int_{\Omega} (\widehat{w} \cdot \nabla)\widehat{w} \cdot \widehat{J} \right| < \eta \, .
$$
Given any vector field $\varphi \in \mathcal{C}^{\infty}_{0}(\Omega)$ (not necessarily divergence-free) and the relative capacity potential $\phi_{\varepsilon} \in H^{2}(\Omega_{\varepsilon})$ (see Lemma \ref{refcapacitypro}), we may use $(1-\phi_{\varepsilon})\varphi \in H^{1}_{0}(\Omega_{\varepsilon})$ as a test function in \eqref{weaksolutionshp} and divide the resulting identity by $\mathcal{Z}^{2}_{\varepsilon}$ to obtain
\begin{equation} \label{bypartsn444}
	\begin{aligned}
	& \dfrac{\eta}{\mathcal{Z}_{\varepsilon}} \int_{\Omega} \nabla \widehat{w_{\varepsilon}} \cdot \nabla \varphi - \dfrac{\eta}{\mathcal{Z}_{\varepsilon}} \int_{\Omega} \nabla \widehat{w_{\varepsilon}} \cdot \nabla (\phi_{\varepsilon} \varphi) + \int_{\Omega} (\widehat{w_{\varepsilon}} \cdot \nabla)\widehat{w_{\varepsilon}} \cdot \varphi - \int_{\Omega} (\widehat{w_{\varepsilon}} \cdot \nabla)\widehat{w_{\varepsilon}} \cdot \phi_{\varepsilon}\varphi  \\[6pt]
	& + \dfrac{1}{\mathcal{Z}_{\varepsilon}} \int_{\Omega} (\widehat{w_{\varepsilon}} \cdot \nabla) \widetilde{J_{\varepsilon}} \cdot (1-\phi_{\varepsilon})\varphi + \dfrac{1}{\mathcal{Z}_{\varepsilon}} \int_{\Omega} (\widetilde{J_{\varepsilon}} \cdot \nabla)\widehat{w_{\varepsilon}} \cdot (1-\phi_{\varepsilon})\varphi - \int_{\Omega} \widehat{p_{\varepsilon}} (\nabla \cdot \varphi) + \int_{\Omega} \widehat{p_{\varepsilon}} (\nabla \cdot \phi_{\varepsilon}\varphi) \\[6pt]
	& \hspace{-4mm} = \dfrac{1}{\mathcal{Z}^{2}_{\varepsilon}} \left( \int_{\Omega} f \cdot (1-\phi_{\varepsilon})\varphi - \eta \int_{\Omega} \nabla \widetilde{J_{\varepsilon}} \cdot \nabla ((1-\phi_{\varepsilon})\varphi) - \int_{\Omega} (\widetilde{J_{\varepsilon}} \cdot \nabla)\widetilde{J_{\varepsilon}} \cdot (1-\phi_{\varepsilon})\varphi \right) \, .
\end{aligned}
\end{equation}
for every $\varepsilon \in I_{*}$, along the subsequences \eqref{convergencesn11}. The convergences in \eqref{convergencesn11} guarantee that
\begin{equation} \label{nsstokeslambda0}
	\lim_{\varepsilon \to 0^{+}} \int_{\Omega} (\widehat{w_{\varepsilon}} \cdot \nabla)\widehat{w_{\varepsilon}} \cdot \varphi = \int_{\Omega} (\widehat{w} \cdot \nabla)\widehat{w} \cdot \varphi \qquad \text{and} \qquad \lim_{\varepsilon \to 0^{+}} \int_{\Omega} \widehat{p_{\varepsilon}} (\nabla \cdot \varphi) = \int_{\Omega} \widehat{p} (\nabla \cdot \varphi) \, .
\end{equation}
On the other hand, applying H\"older and Sobolev inequalities in $\Omega$, from \eqref{cap2} we notice that
\begin{equation} \label{nsstokeslambda2}
	\begin{aligned}
		& \left| \int_{\Omega} \nabla \widehat{w_{\varepsilon}}  \cdot ( \nabla \phi_{\varepsilon} \otimes \varphi + \phi_{\varepsilon} \nabla \varphi) \right| \leq \dfrac{C}{\sigma_{\varepsilon}} \| \nabla \widehat{w_{\varepsilon}} \|_{L^{2}(\Omega)} \| \varphi \|_{W^{1,\infty}(\Omega)} \qquad \forall \varepsilon \in I_{*} \, , \\[6pt]
		& \left| \int_{\Omega} (\widehat{w_{\varepsilon}} \cdot \nabla)\widehat{w_{\varepsilon}} \cdot \phi_{\varepsilon} \varphi \right| \leq \dfrac{C}{\sigma_{\varepsilon}} \| \nabla \widehat{w_{\varepsilon}} \|_{L^{2}(\Omega)} \| \widehat{w_{\varepsilon}} \|_{L^{4}(\Omega)} \| \varphi \|_{L^{\infty}(\Omega)} \qquad \forall \varepsilon \in I_{*} \, , \\[6pt]
& \left| \int_{\Omega} (\widehat{w_{\varepsilon}} \cdot \nabla) \widetilde{J_{\varepsilon}} \cdot (1-\phi_{\varepsilon})\varphi  \right| \leq C \| \nabla \widetilde{J_{\varepsilon}} \|_{L^{2}(\Omega)} \| \widehat{w_{\varepsilon}} \|_{L^{4}(\Omega)} \| \varphi \|_{L^{\infty}(\Omega)} \qquad \forall \varepsilon \in I_{*} \, , \\[6pt]
		& \left| \int_{\Omega} \widehat{p_{\varepsilon}} \left[ \phi_{\varepsilon}(\nabla \cdot \varphi) + \nabla \phi_{\varepsilon} \cdot \varphi  \right] \right| \leq \dfrac{C}{\sigma_{\varepsilon}} \| \widehat{p_{\varepsilon}} \|_{L^{2}(\Omega)} \| \varphi \|_{W^{1,\infty}(\Omega)} \qquad \forall \varepsilon \in I_{*} \, , \\[6pt]
		& \left| \int_{\Omega} f \cdot (1-\phi_{\varepsilon}) \varphi \right| \leq C \| f \|_{L^{2}(\Omega)} \| \varphi \|_{L^{\infty}(\Omega)} \qquad \forall \varepsilon \in I_{*} \, ,
	\end{aligned}
\end{equation}
and the remaining terms appearing in \eqref{bypartsn444} can be estimated in a similar fashion. Observing \eqref{sigma_eps}- \eqref{nsstokeslambda0}-\eqref{nsstokeslambda2}, one can take the limit as $\varepsilon \to 0^{+}$ in \eqref{bypartsn444} to deduce that
$$
\int_{\Omega} (\widehat{w} \cdot \nabla) \widehat{w} \cdot \varphi - \int_{\Omega} \widehat{p} (\nabla \cdot \varphi) = 0 \qquad \forall \varphi \in \mathcal{C}^{\infty}_{0}(\Omega; \mathbb{R}^{2}) \, ,
$$
that is, the pair $(\widehat{w}, \widehat{p}) \in H^{1}_{0}(\Omega) \times L^{2}(\Omega)$ satisfies in distributional form the steady-state Euler equation
\begin{equation} \label{limiteulerlip}
(\widehat{w} \cdot \nabla) \widehat{w} + \nabla \widehat{p} = 0 \, , \quad  \nabla\cdot \widehat{w}=0 \ \ \mbox{ in } \ \ \Omega \, .
\end{equation}
Since $(\widehat{w} \cdot \nabla) \widehat{w} \in L^{3/2}(\Omega)$ (by Sobolev embedding), \eqref{limiteulerlip} proves that actually $\widehat{p} \in W^{1,3/2}(\Omega)$, so that the Euler equation \eqref{limiteulerlip} is satisfied in strong form. Since $\widehat{w} = 0$ on $\partial \Omega$, the Bernoulli law \cite[Lemma 4]{kapitanskii1983spaces} (see \cite[Theorem 2.2]{amick1984existence} and \cite[Theorem 1]{korobkov2011bernoulli} as well) states that $\widehat{p}$ must be constant on each of the connected components of $\partial \Omega$. Thus, there exist $\widehat{p}_{0},...,\widehat{p}_{M} \in \mathbb{R}$ such that
\begin{equation} \label{bernoulli1}
\widehat{p} = \widehat{p}_{j} \quad \text{almost everywhere on} \ \ \partial \Omega_{j} \, , \ \forall j \in \{0,...,M \} \, .
\end{equation}
Given that the scalar pressure $\widehat{p}$ can be defined up to an additive constant, without loss of generality we may assume that $\widehat{p}_{0}=0$. If we multiply the first equation in \eqref{limiteulerlip} by $\widehat{J}$, integrate by parts in $\Omega$ and enforce \eqref{bernoulli1}, the following identity is obtained:
\begin{equation} \label{bernoulli2}
\int_{\Omega} (\widehat{w} \cdot \nabla)\widehat{w} \cdot \widehat{J} = -\int_{\partial \Omega} \widehat{p} (\widehat{J} \cdot \nu ) = - \sum_{j=1}^{M} \widehat{p}_{j} \, F_{j} \, .
\end{equation}
We now proceed as in the proof of \cite[Lemma 4]{neustupa2010new} (see also \cite[Theorem 2.2]{korobkov2014flux}). Define the functions $\widehat{W} \in H_{0}^{1}(\Omega_{0})$ and $\widehat{P} \in W^{1,3/2}(\Omega_{0})$ by
$$
\widehat{W} \doteq
\begin{cases}
	\widehat{w} & \quad \text{in} \ \ \Omega \, ,\\[4pt]
	0 & \quad  \text{in} \ \ \bigcup_{i=1}^{M} \overline{\Omega_{i}}  \, ,
\end{cases}
\qquad \qquad
\widehat{P} \doteq
\begin{cases}
	\widehat{p} & \quad \text{in} \ \ \Omega \, ,\\[4pt]
	\widehat{p}_{j} & \quad  \text{in} \ \ \overline{\Omega_{j}} \quad \forall j \in \{1,...,M \}  \, ,
\end{cases}
$$
so that \eqref{limiteulerlip} implies 
\begin{equation} \label{limiteulerlip2}
	(\widehat{W} \cdot \nabla) \widehat{W} + \nabla \widehat{P} = 0 \, , \quad  \nabla\cdot \widehat{W}=0 \ \ \mbox{ in } \ \ \Omega_{0} \setminus \overline{\Omega_{j}} \, , \ \ \forall j \in \{1,...,M \} \, .
\end{equation}
Given $j \in \{1,...,M \}$, a standard integration by parts and the properties of $\chi_{j}$ in \eqref{harcap1}-\eqref{harcap2} give us
\begin{equation} \label{bernoulli3}
	\int_{\Omega_{0} \setminus \overline{\Omega_{j}}} \nabla \widehat{P} \cdot\nabla \chi_{j} = \int_{\partial \Omega_{j}} \widehat{P} \, \dfrac{\partial \chi_{j}}{\partial \nu} = \widehat{p}_{j} \, \text{Cap}_{\Omega_{0}}(\Omega_{j}) \, ,
\end{equation}
and also
\begin{equation} \label{bernoulli4}
	\int_{\Omega_{0} \setminus \overline{\Omega_{j}}} (\widehat{W} \cdot \nabla) \widehat{W} \cdot \nabla \chi_{j} = - \int_{\Omega_{0} \setminus \overline{\Omega_{j}}} \chi_{j} \left( \nabla \widehat{W} \cdot (\nabla \widehat{W})^{\top} \right) = - \int_{\Omega} \chi_{j} \left( \nabla \widehat{w} \cdot (\nabla \widehat{w})^{\top} \right) \, .
\end{equation}
From \eqref{limiteulerlip2}-\eqref{bernoulli3}-\eqref{bernoulli4} and the Maximum Principle we immediately deduce that
$$
\left| \widehat{p}_{j} \, \text{Cap}_{\Omega_{0}}(\Omega_{j}) \right| = \left| \int_{\Omega} \chi_{j} \left( \nabla \widehat{w} \cdot (\nabla \widehat{w})^{\top} \right) \right| \leq \| \nabla \widehat{w} \|^{2}_{L^{2}(\Omega)} \leq 1 \qquad \forall j \in \{1,...,M \} \, .
$$
Once inserted into \eqref{bernoulli2}, this yields
$$
\left| \int_{\Omega} (\widehat{w} \cdot \nabla)\widehat{w} \cdot \widehat{J} \right| \leq  \sum_{j=1}^{M} \left| \widehat{p}_{j} \, F_{j} \right| \leq \sum_{j=1}^{M} \dfrac{|F_{j} |}{\textup{Cap}_{\Omega_{0}}(\Omega_{j})} < \eta
$$
as a consequence of \eqref{neustupa}, in contradiction with \eqref{contradick}. This concludes the proof.
\end{proof}

Theorem \ref{epslevel} deserves some remarks and comments.

\begin{remark} \label{rem1}
Inequality \eqref{neustupa} imposes a smallness condition on the flux of $v_{*} \in H^{3/2}(\partial \Omega)$ across each of the connected components of $\partial \Omega$, but it does not necessarily enforce a bound on the size of $v_{*}$. As an example we can consider the two-dimensional Taylor-Couette problem for the steady motion of a viscous incompressible fluid in the region between two concentric disks (the inner one at rest and the outer one rotating with arbritarily large and constant angular speed $\omega >0$, see \cite[Chapter II]{landau}).
\end{remark}

\begin{remark}
The hypothesis $f \in H^{1}(\Omega)$ of Theorem \ref{epslevel} is needed in \cite[Theorem 1.1]{korobkov2015solution} to ensure the unrestricted weak solvabilty of problem \eqref{nsstokes0} under the general outflow condition \eqref{gof}. Notice that the inequality \eqref{neustupa} is imposed here to achieve the $\varepsilon$-uniform bound \eqref{uboundpf}, whereas it can also be assumed in order to find a generalized solution to \eqref{nsstokes0} relaxing the regularity constraints on the data of the problem (see \cite[Theorem 1]{neustupa2010new} or \cite[Theorem 2.2]{korobkov2014flux}).
\end{remark}

\begin{remark} \label{rem3}
The uniform bound \eqref{uboundpf} of Theorem \ref{epslevel} can be easily achieved under a smallness assumption on the data (boundary velocity and external force). More precisely, if $\Omega_{\varepsilon}$ is the perforated domain \eqref{perfordomain}, it follows from \cite[Theorem IX.4.1]{galdi2011introduction} the existence of a constant $\delta_{*} > 0$ (depending only on $\Omega$, $\eta$, and $\{ \delta_{0}, \delta_{1} \}$, independent of $\varepsilon \in I_{*}$) such that, whenever
$$
\| v_{*}  \|_{H^{3/2}(\Omega)} + \| f  \|_{H^{1}(\Omega)} < \delta_{*} \, ,
$$
then problem \eqref{nsstokes0} has a unique weak solution $u_{\varepsilon} \in H^{1}_{\sigma}(\Omega_{\varepsilon})$ which, moreover, admits the bound
$$
\| \nabla u_{\varepsilon}  \|_{L^{2}(\Omega_{\varepsilon})} \leq C_{*} \, ,
$$
for some constant $C_{*} > 0$ that depends on $\Omega$, $\eta$, $v_{*}$, $f$ and $\{ \delta_{0}, \delta_{1} \}$ (independent of $\varepsilon \in I_{*}$).
\end{remark}

\subsection{$\varepsilon$-Uniform bounds in the symmetric case} \label{simepslevelsec}
Let $\varepsilon \in I_{*}$ be a fixed parameter. In this subsection we suppose that each of the domains $\Omega_{0}, \Omega_{1},...,\Omega_{M}$ intersects the $x$-axis and is symmetric with respect to the $x$-axis. The \textit{symmetric} version of Lemma \ref{bogtype} can be easily derived as a corollary, and it reads:
\begin{corollary} \label{bogtypesim}
	Let $\Omega_{\varepsilon}$ be the symmetric perforated domain \eqref{perfordomainsim} and $v_{*} \in H^{3/2}(\partial \Omega)$ a symmetric vector field satisfying \eqref{gof}. There exists a symmetric vector field $J_{\varepsilon} \in H^{1}(\Omega_{\varepsilon})$ such that
	\begin{equation} \label{vecjesim}
		\left\{
		\begin{aligned}
			& \nabla \cdot J_{\varepsilon} = 0 \ \ \mbox{in} \ \ \Omega_{\varepsilon} \, ; \qquad J_{\varepsilon} = v_{*} \ \ \mbox{on} \ \ \partial \Omega \, ; \\[3pt]
			& J_{\varepsilon} = 0 \ \ \mbox{on} \ \ \partial K_{\varepsilon} \, ; \qquad \| \nabla J_{\varepsilon} \|_{L^{2}(\Omega_{\varepsilon})} \leq C_{*} \| v_{*} \|_{H^{3/2}(\partial \Omega)} \, ,
		\end{aligned}
		\right.
	\end{equation}
	for some constant $C_{*} > 0$ that depends on $\Omega$ and $\lambda_{0}$, but is independent of $\varepsilon \in I_{*}$. In particular, there exists a symmetric  vector field $\widehat{J} \in H^{1}(\Omega)$ such that
	$$
	\nabla \cdot \widehat{J} = 0 \ \ \mbox{in} \ \ \Omega \qquad \text{and} \qquad \widehat{J} = v_{*} \ \ \mbox{on} \ \ \partial \Omega \, ,
	$$
	and a (not relabeled) subsequence $(\widetilde{J_\varepsilon})_{\varepsilon \in I_{*}} \subset H^{1}(\Omega)$ for which
	\begin{equation} \label{jepsconvsim}
		\widetilde{J_\varepsilon} \rightharpoonup \widehat{J} \ \ \ \text{weakly in} \ H^{1}(\Omega) \qquad \text{and} \qquad \widehat{J_\varepsilon} \to \widehat{J} \ \ \ \text{strongly in} \ L^{4}(\Omega) \ \ \ \text{as} \ \ \varepsilon \to 0^{+} \, .
	\end{equation}
\end{corollary}
\noindent
\begin{proof}
Let $\Psi_{\varepsilon} = (\Psi_{1}^{\varepsilon}, \Psi_{2}^{\varepsilon}) \in H^{1}(\Omega_{\varepsilon})$ be the vector field arising from Lemma \eqref{bogtype}, which satisfies \eqref{vecje} (but is not necessarily symmetric). We then define the vector field $J_{\varepsilon} \in H^{1}(\Omega_{\varepsilon})$ by
$$
J_{\varepsilon}(x,y) \doteq \dfrac{1}{2} ( \Psi_{1}^{\varepsilon}(x,y) + \Psi_{1}^{\varepsilon}(x,-y) , \Psi_{2}^{\varepsilon}(x,y) - \Psi_{2}^{\varepsilon}(x,-y)) \qquad \text{for a.e.} \ (x,y) \in \Omega_{\varepsilon} \, ,
$$
which is symmetric and satisfies all the properties listed in \eqref{vecjesim}. To see this, a direct computation (involving Young's inequality) shows that $\| \nabla J_{\varepsilon} \|_{L^{2}(\Omega_{\varepsilon})} \leq \| \nabla \Psi_{\varepsilon} \|_{L^{2}(\Omega_{\varepsilon})}$. The rest of the proof follows exactly the lines of the proof of Lemma \eqref{bogtype}, where, in addition, we point out that the limit vector field $\widehat{J} \in H^{1}(\Omega)$ is symmetric as a consequence of the strong convergence in \eqref{jepsconvsim}.
\end{proof}

The second main result of this section provides uniform bounds (with respect to $\varepsilon \in I_{*}$) for the symmetric solutions of problem \eqref{nsstokes0}.

\begin{theorem} \label{epslevelsim}
	Let $\Omega_{\varepsilon}$ be the symmetric perforated domain \eqref{perfordomainsim}. For any symmetric vector fields $f \in L^{2}(\Omega)$ and $v_{*} \in H^{3/2}(\partial \Omega)$ such that \eqref{gof} holds, there exists at least one symmetric weak solution $u_{\varepsilon} \in H^{2}(\Omega_{\varepsilon}) \cap H^{1}_{\sigma}(\Omega_{\varepsilon})$ of problem \eqref{nsstokes0}, and an associated symmetric pressure $p_{\varepsilon} \in H^{1}(\Omega_{\varepsilon}) \cap L^{2}_{0}(\Omega_{\varepsilon})$ such that the pair $(u_{\varepsilon},p_{\varepsilon})$ satisfies \eqref{nsstokes0} in strong form. Moreover, the uniform bound
	\begin{equation} \label{uboundpfsim}
		\sup_{\varepsilon \in I_{*}} \left( \| \nabla u_{\varepsilon} \|_{L^{2}(\Omega_{\varepsilon})} + \| p_{\varepsilon} \|_{L^{2}(\Omega_{\varepsilon})} \right) \leq C_{*} \, ,
	\end{equation}
	holds for some constant $C_{*} > 0$ that depends on $\Omega$, $\eta$, $f$, $v_{*}$ and $\lambda_{0}$, but is independent of $\varepsilon \in I_{*}$.
\end{theorem}
\noindent
\begin{proof}
	In what follows, $C > 0$ will always denote a generic constant that depends on $\Omega$, $\eta$ and $\lambda_{0}$ (independently of $\varepsilon \in I_{*}$), but that may change from line to line.
	\par
	Then, given any symmetric vector fields $f \in L^{2}(\Omega)$, $v_{*} \in H^{3/2}(\partial \Omega)$ satisfying \eqref{gof} and $\varepsilon \in I_{*}$, a direct extension of \cite[Theorem 1.1]{amick1984existence} ensures the existence of at least one symmetric weak solution $u_{\varepsilon} \in H^{2}(\Omega_{\varepsilon}) \cap H^{1}_{\sigma}(\Omega_{\varepsilon})$ of problem \eqref{nsstokes0} and an associated symmetric pressure $p_{\varepsilon} \in H^{1}(\Omega_{\varepsilon}) \cap L^{2}_{0}(\Omega_{\varepsilon})$ such that the pair $(u_{\varepsilon},p_{\varepsilon})$ solves \eqref{nsstokes0} in strong form. This follows simply from the fact that the sought solution is required to have zero flux across each of the holes $K^{\varepsilon}_{m}$, for $m \in \mathcal{M}_{\varepsilon}$.  Let $J_{\varepsilon} \in H^{1}(\Omega_{\varepsilon})$ be the symmetric vector field arising from Corollary \eqref{bogtypesim}, which satisfies \eqref{vecjesim}. Set $w_{\varepsilon} \doteq u_{\varepsilon} - J_{\varepsilon} \in H^{1}_{0,\sigma}(\Omega_{\varepsilon})$ (which is a symmetric vector field), so that, in view of the bound in \eqref{vecjesim}$_2$, it suffices to show that
$$
		\sup_{\varepsilon \in I_{*}} \left( \| \nabla w_{\varepsilon} \|_{L^{2}(\Omega_{\varepsilon})} + \| p_{\varepsilon} \|_{L^{2}(\Omega_{\varepsilon})} \right) \leq C \, .
$$
Exactly as in the proof of Theorem \ref{epslevel} we derive identity \eqref{test2} and the estimate \eqref{estpresind}. By contradiction, suppose now that the norms $\| \nabla w_{\varepsilon} \|_{L^{2}(\Omega_{\varepsilon})}$ are not uniformly bounded with respect to $\varepsilon \in I_{*}$. Then, there must exists a subsequence (not being relabeled) such that
\begin{equation} \label{divergentsim}
	\lim_{\varepsilon \to 0^{+}} \mathcal{Z}_{\varepsilon} = + \infty  \quad \text{with} \ \ \mathcal{Z}_{\varepsilon} \doteq \| \nabla w_{\varepsilon} \|_{L^{2}(\Omega_{\varepsilon})} \quad \forall \varepsilon \in I_{*} \, .
\end{equation}
The estimate \eqref{estpresind} enables us to establish that, along this divergent subsequence \eqref{divergentsim}, the following sequences are uniformly bounded with respect to $\varepsilon \in I_{*}$:
$$
( \widehat{w_\varepsilon} )_{\varepsilon \in I_{*}} \doteq \left( \dfrac{\widetilde{w_{\varepsilon}}}{\mathcal{Z}_{\varepsilon}} \right)_{\varepsilon \in I_{*}} \subset H_{0}^{1}(\Omega) \qquad \text{and} \qquad ( \widehat{p_\varepsilon} )_{\varepsilon \in I_{*}} \doteq  \left( \dfrac{\widetilde{p_{\varepsilon}}}{\mathcal{Z}^{2}_{\varepsilon}} \right)_{\varepsilon \in I_{*}} \subset L^{2}(\Omega) \, .
$$
Then, there must exist $\widehat{w} \in H_{0}^{1}(\Omega)$ and $\widehat{p} \in L^{2}(\Omega)$ such that the following convergences hold:
\begin{equation} \label{convergencesn11sim}
	\begin{aligned}
		\widehat{w_\varepsilon} \rightharpoonup \widehat{w} \ \ \ \text{weakly in} \ H^{1}(\Omega) \, ; \qquad &\widehat{w_\varepsilon} \to \widehat{w} \ \ \ \text{strongly in} \ L^{4}(\Omega) \, ; \qquad \widehat{p_\varepsilon} \rightharpoonup \widehat{p} \ \ \ \text{weakly in} \ L^{2}(\Omega) \, ,
	\end{aligned}
\end{equation}
as $\varepsilon \to 0^{+}$, along subsequences that are not being relabeled. Clearly we have $\| \nabla \widehat{w} \|_{L^{2}(\Omega)} \leq 1$, and as in the final part of the proof of Lemma \eqref{bogtype} we deduce that $\nabla \cdot \widehat{w} = 0$ almost everywhere in $\Omega$, that is, $\widehat{w} \in H_{0,\sigma}^{1}(\Omega)$. Moreover, the strong convergence in \eqref{convergencesn11sim} ensures that $\widehat{w}$ is a symmetric vector field in $\Omega$. Imitating the steps of the proof of Theorem \ref{epslevel} we also recover identity \eqref{contradick} and infer that the pair $(\widehat{w}, \widehat{p}) \in H^{1}_{0}(\Omega) \times W^{1,3/2}(\Omega)$ satisfies in strong form the stationary Euler equation \eqref{limiteulerlip} in $\Omega$. Again, the Bernoulli law enforces that $\widehat{p}$ must be constant on each of the connected components of $\partial \Omega$: there exist $\widehat{p}_{0},...,\widehat{p}_{M} \in \mathbb{R}$ such that
\begin{equation} \label{bernoulli1sim}
	\widehat{p} = \widehat{p}_{j} \quad \text{almost everywhere on} \ \ \partial \Omega_{j} \, , \ \forall j \in \{0,...,M \} \, .
\end{equation}
If we multiply the first equation in \eqref{limiteulerlip} by $\widehat{J}$, integrate by parts in $\Omega$ and enforce \eqref{bernoulli1sim}, the following identity is obtained:
\begin{equation} \label{bernoulli2sim}
	\int_{\Omega} (\widehat{w} \cdot \nabla)\widehat{w} \cdot \widehat{J} = -\int_{\partial \Omega} \widehat{p} (\widehat{J} \cdot \nu ) = - \sum_{j=0}^{M} \widehat{p}_{j} \, F_{j} \, .
\end{equation}
A contradiction will be reached after proving that $\widehat{p}_{0} = \widehat{p}_{1} = ...=\widehat{p}_{M}$, as once inserted into \eqref{bernoulli2sim} and combined with \eqref{gof}, this would imply
$$
\int_{\Omega} (\widehat{w} \cdot \nabla)\widehat{w} \cdot \widehat{J} = 0 \, ,
$$
thereby violating \eqref{contradick}. In order to do so, we follow closely the proof of \cite[Theorem 2.3]{amick1984existence} (see also \cite[Theorem 2.1]{korobkov2014flux}). Take $\{ \alpha_{j} \}^{M}_{j=0}, \{ \beta_{j} \}^{M}_{j=0} \subset \mathbb{R}$ such that
$$
\alpha_{j} < \beta_{j} \qquad \text{and} \qquad \partial \Omega_{j} \cap \{(x,y) \in \mathbb{R}^{2} \ | \ y=0 \} = \{ (\alpha_{j},0) \, ,  (\beta_{j},0) \} \qquad \forall j \in \{0,...,M \} \, .
$$
Without loss of generality, we may label the sets $\Omega_{0},...,\Omega_{M}$ in such a way that
$$
\alpha_{0} < \alpha_{1} < \beta_{1} <...< \alpha_{M} < \beta_{M} < \beta_{0} \, .
$$
It suffices to show that $\widehat{p}_{0} = \widehat{p}_{1}$, as the remaining equalities can be proved in analogous way. Since $\Omega_{0}$ and $ \Omega_{1}$ intersect the $x$-axis and they are symmetric with respect to the $x$-axis, there exist $\lambda_{*} > 0$ and two functions $g_{0}, g_{1} \in \mathcal{C}^{2}((-\lambda_{*}, \lambda_{*}); \mathbb{R})$ such that, for every $\lambda \in (0, \lambda_{*})$, the following representations hold:
$$
\begin{aligned}
& \partial \Omega_{0} = \{(x,y) \in \mathbb{R}^{2} \ | \ y \in (-\lambda, \lambda) \, , \ \ x = g_{0}(y) \} \quad \text{near the point} \ (\alpha_{0},0) \, , \\[3pt]
& \partial \Omega_{1} = \{(x,y) \in \mathbb{R}^{2} \ | \ y \in (-\lambda, \lambda) \, , \ \ x = g_{1}(y) \} \quad \text{near the point} \ (\alpha_{1},0) \, ,
\end{aligned}
$$
and such that the following inclusion is observed:
$$
\mathcal{A}_{\lambda} \doteq \{(x,y) \in \mathbb{R}^{2} \ | \ 0 < y < \lambda \, , \ \ g_{0}(y) < x < g_{1}(y) \} \subset \Omega \, .
$$
Write $\widehat{w} = (\widehat{w}_{1}, \widehat{w}_{2})$ in $\Omega$ and define the \textit{total-head pressure} as
$$
\widehat{\Phi} \doteq \widehat{p} + \dfrac{1}{2} | \widehat{w} |^2  \quad \text{in} \ \ \Omega \, ,
$$
so that $\widehat{\Phi} \in W^{1,3/2}(\Omega)$ and it clearly verifies the identity
$$
\dfrac{\partial \widehat{\Phi}}{\partial x} = \widehat{w}_{2} \left( \dfrac{\partial \widehat{w}_{2}}{\partial x} - \dfrac{\partial \widehat{w}_{1}}{\partial y}  \right) \quad \text{in} \ \ \Omega \, .
$$
Integrating this last identity over $\mathcal{A}_{\lambda}$ and enforcing \eqref{bernoulli1sim}, we obtain
$$
\lambda (\widehat{p}_{1} - \widehat{p}_{0}) = \int_{\mathcal{A}_{\lambda}} \widehat{w}_{2} \left( \dfrac{\partial \widehat{w}_{2}}{\partial x} - \dfrac{\partial \widehat{w}_{1}}{\partial y}  \right) \, ,
$$
so that
\begin{equation} \label{amick5}
| \widehat{p}_{1} - \widehat{p}_{0} | = \dfrac{1}{\lambda} \left| \int_{\mathcal{A}_{\lambda}} \widehat{w}_{2} \left( \dfrac{\partial \widehat{w}_{2}}{\partial x} - \dfrac{\partial \widehat{w}_{1}}{\partial y}  \right) \right| \qquad \forall \lambda \in (0, \lambda_{*}) \, .
\end{equation}
If we extend $\widehat{w}_{2}$ by zero to the whole set $\mathbb{R} \times (0,\lambda)$, by symmetry we have $\widehat{w}_{2}(x,0)=0$ (in the sense of traces) for every $x \in \mathbb{R}$, so that Hardy's inequality can be applied to yield
\begin{equation} \label{amick6}
\int_{\mathcal{A}_{\lambda}} \dfrac{| \widehat{w}_{2}(x,y) |^2}{y^2} \, dx \, dy = \int\limits^{+\infty}_{-\infty} \int\limits^{\lambda}_{0} \dfrac{| \widehat{w}_{2}(x,y) |^2}{y^2} \, dy \, dx \leq 4 \int\limits^{+\infty}_{-\infty} \int\limits^{\lambda}_{0} \left| \dfrac{\partial \widehat{w}_{2}}{\partial y}(x,y) \right|^{2} dy \, dx \leq 4 \int_{\mathcal{A}_{\lambda}} | \nabla \widehat{w} |^2 \, .
\end{equation}
Combining \eqref{amick5}-\eqref{amick6} (together with H\"older's inequality) leads to the estimate
$$
\begin{aligned}
| \widehat{p}_{1} - \widehat{p}_{0} |^{2} \leq \dfrac{1}{\lambda^2} \left( \int_{\mathcal{A}_{\lambda}} y^2 \dfrac{| \widehat{w}_{2}(x,y) |^2}{y^2} \, dx \, dy \right) \left( \int_{\mathcal{A}_{\lambda}} \left| \dfrac{\partial \widehat{w}_{2}}{\partial x} - \dfrac{\partial \widehat{w}_{1}}{\partial y}  \right|^2  \right) \leq 16 \int_{\mathcal{A}_{\lambda}} | \nabla \widehat{w} |^2 \to 0 \quad \text{as} \ \ \lambda \to 0^{+} \, ,
\end{aligned}
$$
thereby proving that $\widehat{p}_{1} = \widehat{p}_{0}$. This concludes the proof.
\end{proof}

\begin{remark} \label{rem33}
The uniform bound \eqref{uboundpfsim} of Theorem \ref{epslevelsim} can be easily achieved under a smallness assumption on the data (boundary velocity and external force). More precisely, if $\Omega_{\varepsilon}$ is the symmetric perforated domain \eqref{perfordomainsim}, it follows from \cite[Theorem IX.4.1]{galdi2011introduction} the existence of a constant $\delta_{*} > 0$ (depending only on $\Omega$, $\eta$, and $\lambda_{0}$, independent of $\varepsilon \in I_{*}$) such that, whenever
$$
\| v_{*}  \|_{H^{3/2}(\Omega)} + \| f  \|_{L^{2}(\Omega)} < \delta_{*} \, ,
$$
then \eqref{nsstokes0} has a unique (symmetric) weak solution $u_{\varepsilon} \in H^{1}_{\sigma}(\Omega_{\varepsilon})$ which, moreover, admits the bound
$$
\| \nabla u_{\varepsilon}  \|_{L^{2}(\Omega_{\varepsilon})} \leq C_{*} \, ,
$$
for some constant $C_{*} > 0$ that depends on $\Omega$, $\eta$, $v_{*}$, $f$ and $\lambda_{0}$ (independent of $\varepsilon \in I_{*}$).
\end{remark}

\newpage
\section{Asymptotic behavior as $\varepsilon \to 0^{+}$: homogenized equations} \label{energymethod}
By employing the renowned \textit{energy method} of Tartar \cite[Appendix]{sanchez1980non} (see also \cite[Chapter 15]{tartar2009general}), in this section we derive the effective (or \textit{homogenized}) equations satisfied by the solutions of \eqref{nsstokes0} as $\varepsilon \to 0^{+}$.

\begin{theorem} \label{effectiveq1}
	Let $(\Omega_{\varepsilon})_{\varepsilon \in I_{*}}$ be the family of perforated domains \eqref{perfordomain}. Given any $f \in H^{1}(\Omega)$ and $v_{*} \in H^{3/2}(\partial \Omega)$ satisfying \eqref{gof}, let $(u_{\varepsilon},p_{\varepsilon}) \in H_{\sigma}^{1}(\Omega_{\varepsilon}) \times L_{0}^{2}(\Omega_{\varepsilon})$ be a weak solution of  \eqref{nsstokes0}. Then, assuming condition \eqref{neustupa}, up to the extraction of a subsequence, the sequence of extended functions $\{(\widetilde{u_{\varepsilon}},\widetilde{p_{\varepsilon}}) \}_{\varepsilon \in I_{*}} \subset H^{1}(\Omega) \times L_{0}^{2}(\Omega)$ converges strongly to a weak solution $(u,p) \in H^{1}(\Omega) \times L_{0}^{2}(\Omega)$ of problem \eqref{nsstokes0} in $\Omega$ as $\varepsilon \to 0^{+}$. Furthermore, $(u,p) \in H^{2}(\Omega) \times H^{1}(\Omega)$ and it satisfies in strong form the system
\begin{equation} \label{nsstokesfinal}
	\left\{
	\begin{aligned}
		& -\eta\Delta u+(u\cdot\nabla)u+\nabla p=f \, , \quad  \nabla\cdot u=0 \ \ \mbox{ in } \ \ \Omega \, , \\[4pt]
		& u=v_{*} \ \ \mbox{ on } \ \ \partial \Omega \, .
	\end{aligned}
	\right.
\end{equation}
\end{theorem}
\noindent
\begin{proof}
	In what follows, $C > 0$ will always denote a generic constant that depends on $\Omega$, $\eta$ and $\{ \delta_{0}, \delta_{1} \}$ (independently of $\varepsilon \in I_{*}$), but that may change from line to line.
\par
Given any $f \in H^{1}(\Omega)$, $v_{*} \in H^{3/2}(\partial \Omega)$ satisfying \eqref{gof} and $\varepsilon \in I_{*}$, let $(u_{\varepsilon},p_{\varepsilon}) \in H_{\sigma}^{1}(\Omega_{\varepsilon}) \times L_{0}^{2}(\Omega_{\varepsilon})$ be a weak solution of  \eqref{nsstokes0}. From Theorem \ref{epslevel} we know that $(u_{\varepsilon},p_{\varepsilon}) \in H^{2}(\Omega_{\varepsilon}) \times H^{1}(\Omega_{\varepsilon})$ solves \eqref{nsstokes0} in strong form and that $\{(\widetilde{u_{\varepsilon}},{\color{black}\widetilde{p_{\varepsilon}}}) \}_{\varepsilon \in I_{*}} \subset H^{1}(\Omega) \times L_{0}^{2}(\Omega)$. In particular, the weak formulation of Definition \eqref{weaksolution} incorporating the pressure term reads:
\begin{equation} \label{weaksolutionshp2}
\eta \int_{\Omega_{\varepsilon}} \nabla u_{\varepsilon} \cdot \nabla \varphi + \int_{\Omega_{\varepsilon}} (u_{\varepsilon} \cdot \nabla)u_{\varepsilon} \cdot \varphi - \int_{\Omega_{\varepsilon}} p_{\varepsilon} (\nabla \cdot \varphi)  = \int_{\Omega_{\varepsilon}} f \cdot \varphi \qquad \forall \varphi \in H^{1}_{0}(\Omega_{\varepsilon}) \, .
\end{equation}
Now, given any scalar function $\phi \in \mathcal{C}_{0}^{\infty}(\Omega)$, an integration by parts and the divergence-free condition in \eqref{nsstokes0}$_1$ imply that
$$
- \int_{\Omega} \phi (\nabla \cdot \widetilde{u_{\varepsilon}} )  = \int_{\Omega} \widetilde{u_{\varepsilon}} \cdot \nabla \phi = \int_{\Omega_{\varepsilon}} u_{\varepsilon} \cdot \nabla \phi = \int_{\partial \Omega_{\varepsilon}} \phi (u_{\varepsilon} \cdot \nu) = 0 \qquad \forall \varepsilon \in I_{*} \, ,
$$
since $u_{\varepsilon}$ vanishes on $\partial K_{\varepsilon}$ and so does $\phi$ on $\partial \Omega$. This proves $\nabla \cdot \widetilde{u_{\varepsilon}} = 0$ almost everywhere in $\Omega$, that is, $\widetilde{u_{\varepsilon}} \in H_{\sigma}^{1}(\Omega)$ for every $\varepsilon \in I_{*}$. Moreover, \eqref{uboundpf}-\eqref{extvel} ensure that the sequences $(\widetilde{u_{\varepsilon}})_{\varepsilon \in I_{*}} \subset H_{\sigma}^{1}(\Omega)$ and $(\widetilde{p_{\varepsilon}})_{\varepsilon \in I_{*}} \subset L_{0}^{2}(\Omega)$ are uniformly bounded, so there exist $u \in H_{\sigma}^{1}(\Omega)$ and $p \in L_{0}^{2}(\Omega)$ such that the following convergences hold as $\varepsilon \to 0^{+}$:
\begin{equation} \label{convergencesn22}
\begin{aligned}
& \widetilde{u_{\varepsilon}} \rightharpoonup u \ \ \ \text{weakly in} \ H^{1}(\Omega) \, ; \qquad \widetilde{u_{\varepsilon}} \to u \ \ \ \text{strongly in} \ L^{4}(\Omega) \, ; \\[3pt]
& \widetilde{u_{\varepsilon}} \to u \ \ \ \text{strongly in} \ L^{2}(\partial \Omega) \, ; \qquad \widetilde{p_\varepsilon} \rightharpoonup p \ \ \ \text{weakly in} \ L^{2}(\Omega) \, ,
\end{aligned}
\end{equation}
along subsequences that are not being relabeled. The strong convergence in \eqref{convergencesn22}$_2$ ensures that $u = v_{*}$ almost everywhere on $\partial \Omega$. Given any vector field $\varphi \in \mathcal{C}_{0}^{\infty}(\Omega)$ (not necessarily divergence-free) and the relative capacity potential $\phi_{\varepsilon} \in H^{2}(\Omega_{\varepsilon})$ (see Lemma \ref{refcapacitypro}), we may use $(1-\phi_{\varepsilon})\varphi \in H^{1}_{0}(\Omega_{\varepsilon})$ as a test function in \eqref{weaksolutionshp2} to obtain
\begin{equation} \label{tartar1}
	\begin{aligned}
		& \eta \int_{\Omega} \nabla \widetilde{u_{\varepsilon}} \cdot \nabla \varphi - \eta \int_{\Omega} \nabla \widetilde{u_{\varepsilon}} \cdot \nabla (\phi_{\varepsilon} \varphi) + \int_{\Omega} (\widetilde{u_{\varepsilon}} \cdot \nabla)\widetilde{u_{\varepsilon}} \cdot \varphi - \int_{\Omega} (\widetilde{u_{\varepsilon}} \cdot \nabla)\widetilde{u_{\varepsilon}} \cdot \phi_{\varepsilon}\varphi \\[6pt]
		& - \int_{\Omega} \widetilde{p_{\varepsilon}} (\nabla \cdot \varphi) + \int_{\Omega} \widetilde{p_{\varepsilon}} (\nabla \cdot \phi_{\varepsilon}\varphi) =  \int_{\Omega} f \cdot (1-\phi_{\varepsilon}) \varphi \, ,
	\end{aligned}
\end{equation}
for every $\varepsilon \in I_{*}$, along the subsequences \eqref{convergencesn22}. Observing \eqref{cap2}, notice that
\begin{equation} \label{bypartsn6}
	\left| \int_{\Omega} f \cdot (1-\phi_{\varepsilon})\varphi - \int_{\Omega} f \cdot \varphi \right| \leq \left| \int_{K_{\varepsilon}} f \cdot \varphi \right| + \dfrac{C}{\sigma_{\varepsilon}} \| f \|_{L^{2}(\Omega)}  \| \varphi \|_{L^{\infty}(\Omega)} \to 0 \quad \text{as} \ \ \varepsilon \to 0^{+} \, ,
\end{equation}
because $f \cdot \varphi \in L^{1}(\Omega)$ and $| K_{\varepsilon} | \to 0$ as $\varepsilon \to 0^{+}$. Also, with the help of the convergences in \eqref{convergencesn22} (see again \eqref{nsstokeslambda0}-\eqref{nsstokeslambda2}) we can easily prove that
\begin{equation} \label{bypartsn5}
	\begin{aligned}
		& \lim_{\varepsilon \to 0^{+}} \int_{\Omega} \nabla \widetilde{u_{\varepsilon}} \cdot \nabla \varphi = \int_{\Omega} \nabla u \cdot \nabla \varphi \, , \qquad  \lim_{\varepsilon \to 0^{+}} \int_{\Omega} (\widetilde{u_{\varepsilon}} \cdot \nabla)\widetilde{u_{\varepsilon}} \cdot \varphi = \int_{\Omega} (u \cdot \nabla)u \cdot \varphi \, , \\[6pt]
		& \lim_{\varepsilon \to 0^{+}} \int_{\Omega} \nabla \widetilde{u_{\varepsilon}}  \cdot ( \nabla \phi_{\varepsilon} \otimes \varphi + \phi_{\varepsilon} \nabla \varphi) = \lim_{\varepsilon \to 0^{+}} \int_{\Omega} (\widetilde{u_{\varepsilon}} \cdot \nabla)\widetilde{u_{\varepsilon}} \cdot \phi_{\varepsilon}\varphi = 0 \, , \\[6pt]
		& \lim_{\varepsilon \to 0^{+}} \int_{\Omega} \widetilde{p_{\varepsilon}} (\nabla \cdot \varphi) = \int_{\Omega} p (\nabla \cdot \varphi) \, , \qquad \lim_{\varepsilon \to 0^{+}} \int_{\Omega} \widetilde{p_{\varepsilon}} \left[ \phi_{\varepsilon}(\nabla \cdot \varphi) + \nabla \phi_{\varepsilon} \cdot \varphi  \right] = 0 \, .
	\end{aligned}
\end{equation}
In view of \eqref{bypartsn6}-\eqref{bypartsn5}, we then take the limit as $\varepsilon \to 0^{+}$ in \eqref{tartar1} to deduce
$$
\eta \int_{\Omega} \nabla u \cdot \nabla \varphi + \int_{\Omega} (u \cdot \nabla)u \cdot \varphi - \int_{\Omega} p (\nabla \cdot \varphi)  = \int_{\Omega} f \cdot \varphi \qquad \forall \varphi \in \mathcal{C}_{0}^{\infty}(\Omega;\mathbb{R}^{2}) \, ,
$$
so that, by density, $(u,p) \in H_{\sigma}^{1}(\Omega) \times L_{0}^{2}(\Omega)$ is a weak solution of \eqref{nsstokes0} in $\Omega$ according to Definition \ref{weaksolution}. Since $\Omega$ has a boundary of class $\mathcal{C}^2$, $f \in L^{2}(\Omega)$ and $v_{*} \in H^{3/2}(\partial \Omega)$, the usual regularity results for the steady-state Navier-Stokes equations with non-homogeneous Dirichlet boundary conditions (see, for example, \cite[Theorem IX.5.2]{galdi2011introduction}) then ensure that $(u,p) \in H^{2}(\Omega) \times H^{1}(\Omega)$ satisfies in strong form the system \eqref{nsstokesfinal}.
\par
In order to show the strong convergence in $H^{1}(\Omega) \times L^{2}(\Omega)$, in view of the weak convergences in \eqref{convergencesn22}, it clearly suffices to prove that
\begin{equation} \label{strong}
	\lim_{\varepsilon \to 0^{+}} \| \nabla \widetilde{u_{\varepsilon}} \|_{L^{2}(\Omega)} = \| \nabla u \|_{L^{2}(\Omega)} \qquad \text{and} \qquad \lim_{\varepsilon \to 0^{+}} \| \widetilde{p_{\varepsilon}} \|_{L^{2}(\Omega)} = \| p \|_{L^{2}(\Omega)} \, .
\end{equation}
Let $\widehat{J} \in H^{1}(\Omega)$ be the vector field arising from Lemma \ref{bogtype}. We take $u_{\varepsilon} - (1-\phi_{\varepsilon})\widehat{J} \in H^{1}_{0}(\Omega_{\varepsilon})$ as a test function in \eqref{weaksolutionshp2} to get
\begin{equation} \label{tartar2}
	\begin{aligned}
		& \eta \| \nabla \widetilde{u_{\varepsilon}} \|^{2}_{L^{2}(\Omega)} - \eta \int_{\Omega} \nabla \widetilde{u_{\varepsilon}} \cdot \nabla \widehat{J} +  \eta \int_{\Omega} \nabla \widetilde{u_{\varepsilon}} \cdot \nabla (\phi_{\varepsilon} \widehat{J}) + \int_{\Omega_{\varepsilon}} (u_{\varepsilon} \cdot \nabla)u_{\varepsilon} \cdot u_{\varepsilon} - \int_{\Omega} (\widetilde{u_{\varepsilon}} \cdot \nabla)\widetilde{u_{\varepsilon}} \cdot \widehat{J} \\[6pt]
		& + \int_{\Omega} (\widetilde{u_{\varepsilon}} \cdot \nabla)\widetilde{u_{\varepsilon}} \cdot \phi_{\varepsilon} \widehat{J} - \int_{\Omega} \widetilde{p_{\varepsilon}} (\nabla \cdot \phi_{\varepsilon}\widehat{J}) = \int_{\Omega} f \cdot \widetilde{u_{\varepsilon}} -  \int_{\Omega} f \cdot  (1-\phi_{\varepsilon})\widehat{J} \, ,
	\end{aligned}
\end{equation}
for every $\varepsilon \in I_{*}$, along the subsequences \eqref{convergencesn22}. A standard integration by parts yields
\begin{equation} \label{tartar3}
\int_{\Omega_{\varepsilon}} (u_{\varepsilon} \cdot \nabla)u_{\varepsilon} \cdot u_{\varepsilon} = \dfrac{1}{2} \int_{\partial \Omega} |v_{*}|^{2} (v_{*} \cdot \nu) \qquad \forall \varepsilon \in I_{*} \, .
\end{equation}
Upon substitution of \eqref{tartar3} into \eqref{tartar2}, and arguing as in \eqref{bypartsn6}-\eqref{bypartsn5}, we can take the limit as $\varepsilon \to 0^{+}$ in \eqref{tartar2} to infer the identity
\begin{equation} \label{tartar4}
\eta \lim_{\varepsilon \to 0^{+}} \| \nabla \widetilde{u_{\varepsilon}} \|^{2}_{L^{2}(\Omega)} = \eta \int_{\Omega} \nabla u \cdot \nabla \widehat{J} - \dfrac{1}{2} \int_{\partial \Omega} |v_{*}|^{2} (v_{*} \cdot \nu) + \int_{\Omega} (u \cdot \nabla)u \cdot \widehat{J} + \int_{\Omega} f \cdot (u -  \widehat{J}) \, .
\end{equation}
On the other hand, multiplying the first identity in \eqref{nsstokesfinal}$_1$ by $u - \widehat{J} \in H^{1}_{0}(\Omega)$ and integrating by parts in $\Omega$ (arguing again as in \eqref{tartar3}), entails
\begin{equation} \label{tartar5}
\eta \| \nabla u \|^{2}_{L^{2}(\Omega)} = \eta \int_{\Omega} \nabla u \cdot \nabla \widehat{J} - \dfrac{1}{2} \int_{\partial \Omega} |v_{*}|^{2} (v_{*} \cdot \nu) + \int_{\Omega} (u \cdot \nabla)u \cdot \widehat{J} + \int_{\Omega} f \cdot (u -  \widehat{J}) \, ,
\end{equation}
so that \eqref{tartar4}-\eqref{tartar5} secure the first equality in \eqref{strong}, that is, $\widetilde{u_{\varepsilon}} \to u$ strongly in $H^{1}(\Omega)$ as $\varepsilon \to 0^{+}$. Now, since $p_{\varepsilon} \in L^{2}_{0}(\Omega_{\varepsilon})$, let $P_{\varepsilon} \in H_{0}^{1}(\Omega_{\varepsilon})$ be a vector field verifying \eqref{vecjepre}. Observing that $\widetilde{P_{\varepsilon}} \in H_{0}^{1}(\Omega)$ and that $\| \nabla \widetilde{P_{\varepsilon}} \|_{L^{2}(\Omega)}=\| \nabla P_{\varepsilon} \|_{L^{2}(\Omega_{\varepsilon})}$, \eqref{uboundpf} then ensures that the sequence $(\widetilde{P_{\varepsilon}})_{\varepsilon \in I_{*}} \subset H_{0}^{1}(\Omega)$ is uniformly bounded, and so we deduce the existence of $P \in H_{0}^{1}(\Omega)$ such that
\begin{equation} \label{convergencesj}
\widetilde{P_\varepsilon} \rightharpoonup P \ \ \ \text{weakly in} \ H^{1}(\Omega) \qquad \text{and} \qquad \widetilde{P_\varepsilon} \to P \ \ \ \text{strongly in} \ L^{4}(\Omega) \quad \text{as} \ \ \varepsilon \to 0^{+} \, ,
\end{equation}
along a (not relabeled) subsequence. Given any scalar function $\phi \in \mathcal{C}_{0}^{\infty}(\Omega)$, since $\nabla \cdot P_{\varepsilon} = p_{\varepsilon}$ in $\Omega_{\varepsilon}$ and $P_{\varepsilon}$ vanishes on $\partial \Omega_{\varepsilon}$ for any $\varepsilon \in I_{*}$, an integration by parts gives us
$$
\int_{\Omega} \widetilde{P_{\varepsilon}} \cdot \nabla \phi = \int_{\Omega_{\varepsilon}} P_{\varepsilon} \cdot \nabla \phi = - \int_{\Omega_{\varepsilon}} \phi \, p_{\varepsilon} = - \int_{\Omega} \phi \, \widetilde{p_{\varepsilon}} \qquad \forall \varepsilon \in I_{*} \, ,
$$
along the subsequences \eqref{convergencesn22}-\eqref{convergencesj}. Taking the limit in this last identity as $\varepsilon \to 0^{+}$ entails
$$
\int_{\Omega} P \cdot \nabla \phi = - \int_{\Omega} \phi \, p \qquad \forall \phi \in \mathcal{C}_{0}^{\infty}(\Omega; \mathbb{R}) \, ,
$$
that is, $\nabla \cdot P = p$ in $\Omega$. As in the proof of Theorem \ref{epslevel} we derive the identity
\begin{equation} \label{tartar6}
\| p_{\varepsilon} \|^{2}_{L^{2}(\Omega_{\varepsilon})} = \eta \int_{\Omega} \nabla \widetilde{u_{\varepsilon}} \cdot \nabla \widetilde{P_{\varepsilon}} + \int_{\Omega} (\widetilde{u_{\varepsilon}} \cdot \nabla)\widetilde{u_{\varepsilon}} \cdot \widetilde{P_{\varepsilon}} - \int_{\Omega} f \cdot \widetilde{P_{\varepsilon}} \qquad \forall \varepsilon \in I_{*} \, ,
\end{equation}
along the subsequences \eqref{convergencesn22}-\eqref{convergencesj}. Similarly, multiply the first identity in \eqref{nsstokesfinal}$_1$ by $P$ and integrate by parts in $\Omega$ to provide
$$
\| p \|^{2}_{L^{2}(\Omega)} = \eta \int_{\Omega} \nabla u \cdot \nabla P + \int_{\Omega} (u \cdot \nabla)u \cdot P - \int_{\Omega} f \cdot P \, .
$$
Recalling that $\widetilde{u_{\varepsilon}} \to u$ strongly in $H^{1}(\Omega)$ as $\varepsilon \to 0^{+}$, we can take the limit in \eqref{tartar6} as $\varepsilon \to 0^{+}$ to reach the second equality in \eqref{strong}. This concludes the proof.
\end{proof}

Concerning the symmetric case, we have the following result, analogous to Theorem \ref{effectiveq1}, whose proof is omitted (for the sake of brevity, since it is very similar to the proof of Theorem  \ref{effectiveq1} with obvious minor modifications):

\begin{theorem} \label{effectiveq2}
Let $(\Omega_{\varepsilon})_{\varepsilon \in I_{*}}$ be the family of symmetric perforated domains \eqref{perfordomainsim}. For any given symmetric vector fields $f \in H^{1}(\Omega)$ and $v_{*} \in H^{3/2}(\partial \Omega)$ satisfying \eqref{gof}, let $(u_{\varepsilon},p_{\varepsilon}) \in H_{\sigma}^{1}(\Omega_{\varepsilon}) \times L_{0}^{2}(\Omega_{\varepsilon})$ be a symmetric weak solution of \eqref{nsstokes0}. Then, up to the extraction of a subsequence, the sequence of extended functions $\{(\widetilde{u_{\varepsilon}},\widetilde{p_{\varepsilon}}) \}_{\varepsilon \in I_{*}} \subset H^{1}(\Omega) \times L_{0}^{2}(\Omega)$ converges strongly to a symmetric weak solution $(u,p) \in H^{1}(\Omega) \times L_{0}^{2}(\Omega)$ of problem \eqref{nsstokes0} in $\Omega$ as $\varepsilon \to 0^{+}$. Furthermore, $(u,p) \in H^{2}(\Omega) \times H^{1}(\Omega)$ and it solves in strong form the system
$$
	\left\{
	\begin{aligned}
		& -\eta\Delta u+(u\cdot\nabla)u+\nabla p=f \, , \quad  \nabla\cdot u=0 \ \ \mbox{ in } \ \ \Omega \, , \\[4pt]
		& u=v_{*} \ \ \mbox{ on } \ \ \partial \Omega \, .
	\end{aligned}
	\right.
$$
\end{theorem}

\newpage
\section{Remarks on the fully general two-dimensional case} \label{remarksfulcase}
Let $\varepsilon \in I_{*}$ and $\Omega_{\varepsilon}$ be as in \eqref{perfordomain}, where we additionally suppose that the sequence of sets $( K_{n} )_{n \in \mathbb{N}}$ is uniformly Lipschitz ($\bigstar$). Given any $f \in H^{1}(\Omega)$ and $v_{*} \in H^{3/2}(\partial \Omega)$ satisfying \eqref{gof}, \cite[Theorem 1.1]{korobkov2015solution} (see also \cite[Remark 1.1]{korobkov2015solution}) ensures the existence of at least one weak solution $u_{\varepsilon} \in H^{2}(\Omega_{\varepsilon}) \cap H^{1}_{\sigma}(\Omega_{\varepsilon})$ of problem \eqref{nsstokes0} and an associated pressure $p_{\varepsilon} \in H^{1}(\Omega_{\varepsilon}) \cap L^{2}_{0}(\Omega_{\varepsilon})$ such that the pair $(u_{\varepsilon}, p_{\varepsilon})$ solves \eqref{nsstokes0} in strong form. In order to prove this result (which does not impose any additional assumption on the data of the problem or on the size of the boundary fluxes), the authors of \cite{korobkov2015solution} use the Leray-Schauder Principle, where the required uniform bounds (with respect to some parameter $\lambda \in [0,1]$) are reached through a contradiction argument that employs Bernoulli's law for weak solutions of the stationary Euler equations and a generalization of the Morse-Sard Theorem for Sobolev functions. An essential ingredient in the proof of \cite[Theorem 1.1]{korobkov2015solution} is the uniform boundedness (with respect to $\lambda \in [0,1]$) of the $W^{1,3/2}(\Omega_{\varepsilon})$-norm of the scalar pressure, see \cite[Lemma 3.1]{korobkov2015solution}. This would indicate that, in order to obtain $\varepsilon$-uniform bounds for the solutions of \eqref{nsstokes0} following the method given in \cite[Theorem 1.1]{korobkov2015solution}, one would need to build a $W^{1,3/2}(\Omega)$-uniform extension for the pressure $p_{\varepsilon} \in H^{1}(\Omega_{\varepsilon}) \cap L^{2}_{0}(\Omega_{\varepsilon})$ inside the holes $K_{\varepsilon}$. Inspired by \cite[Lemma 3.1]{juodagalvyte2020time} and \cite[Lemma 1.7]{masmoudi2002homogenization}, the purpose of this section is precisely to illustrate the fact that such a uniform extension cannot be achieved in a simple way. A first explanation derives from the observation that, even though $p_{\varepsilon} \in H^{1}(\Omega_{\varepsilon})$, the system \eqref{nsstokes0} does not provide any information concerning the behavior of $p_{\varepsilon}$ on the boundary of the holes $\partial K_{\varepsilon}$. A second, more convincing, explanation involves a \textit{microscopic} analysis of the boundary-value problem \eqref{nsstokes0} near each single hole $K^{\varepsilon}_{n}$, with $n \in \{1,...,N(\varepsilon) \}$. We will show that
\begin{equation} \label{masmoudi0}
 \| \nabla p_{\varepsilon} \|_{L^{3/2}(\Omega_{\varepsilon})} \leq \dfrac{C_{*}}{a_{\varepsilon}} \left(1 + \| f \|^{2}_{L^{2}(\Omega)} + \| \nabla u_{\varepsilon} \|^{4}_{L^{2}(\Omega_{\varepsilon})} + \| p_{\varepsilon} \|^{2}_{L^{2}(\Omega_{\varepsilon})} + \| v_{*} \|^{4}_{H^{3/2}(\partial \Omega)} \right) \qquad \forall \varepsilon \in I_{*} \, ,
\end{equation}
for some constant $C_{*}>0$ independent of $\varepsilon \in I_{*}$. To do so, after translation and rescaling, \eqref{perforation}$_1$ implies
$$
\overline{K_{n}} \subset D((0,0), \delta_{0} ) \subset D((0,0), \delta_{1} ) \subset D \left( (0,0), \dfrac{\delta_{1} \varepsilon}{a_{\varepsilon}} \right) \qquad \forall n \in \{1,...,N(\varepsilon)\} \, ,
$$
since $a_{\varepsilon} \ll \varepsilon$ for every $\varepsilon \in I_{*}$. Notice then that
\begin{equation}\label{cv}
z \in  D \left(\xi^{\varepsilon}_{n}, \delta_{0} a_{\varepsilon} \right) \setminus \overline{K_{n}^{\varepsilon}} \quad \Longleftrightarrow \quad \dfrac{z-\xi^{\varepsilon}_{n}}{a_{\varepsilon}} \in D((0,0), \delta_{0} ) \setminus \overline{K_{n}} \qquad \forall n \in \{1,...,N(\varepsilon)\} \, .
\end{equation}
Fix any $n \in \{1,...,N(\varepsilon)\}$. In what follows, $C > 0$ will always denote a generic constant that is independent of $\varepsilon \in I_{*}$ and $n \in \{1,...,N(\varepsilon)\}$, but that may change from line to line. Set $D^n_{j} \doteq D((0,0), \delta_{j} ) \setminus \overline{K_{n}}$, for $j \in \{0,1\}$, and define $(U_{\varepsilon}, P_{\varepsilon},F) \in H^{2}(D^n_1) \times H^{1}(D^n_1) \times H^{1}(D^n_1)$ by
$$
U_{\varepsilon}(z) \doteq u_{\varepsilon}(a_{\varepsilon} z + \xi^{\varepsilon}_{n}) \, , \quad P_{\varepsilon}(z) \doteq p_{\varepsilon}(a_{\varepsilon} z + \xi^{\varepsilon}_{n}) \, , \quad F(z) \doteq f(a_{\varepsilon} z + \xi^{\varepsilon}_{n}) \qquad \forall z \in D^n_1 \, ,
$$
so that these functions satisfy the following Stokes-type system in $D^n_1$:
$$
	\left\{
	\begin{aligned}
		& -\eta \Delta U_{\varepsilon} + a_{\varepsilon} \nabla P_{\varepsilon}= a^{2}_{\varepsilon} F - a_{\varepsilon} (U_{\varepsilon}\cdot\nabla)U_{\varepsilon} \, , \quad  \nabla\cdot U_{\varepsilon} = 0 \ \ \mbox{ in } \ \ D^n_1 \, , \\[4pt]
		& U_{\varepsilon} = 0 \ \ \mbox{ on } \ \ \partial K_{n} \, .
	\end{aligned}
	\right.
$$
The usual local regularity estimates for the Stokes equations (see \cite[Teorema, page 311]{cattabriga1961problema} or \cite[Theorem IV.4.1]{galdi2011introduction}) entail
\begin{equation} \label{regularity1}
\| \nabla P_{\varepsilon} \|^{3/2}_{L^{3/2}(D^n_0)} \leq C \left( \| a^{2}_{\varepsilon} F - a_{\varepsilon} (U_{\varepsilon}\cdot\nabla)U_{\varepsilon} \|^{3/2}_{L^{3/2}(D^n_1)} + \| \nabla U_{\varepsilon} \|^{3/2}_{L^{3/2}(D^n_1)} + \| P_{\varepsilon} \|^{3/2}_{L^{3/2}(D^n_1)} \right) \, ,
\end{equation}
where we emphasize that, as a consequence of property ($\bigstar$), the constant $C>0$ entering \eqref{regularity1} can be bounded independently of $n \in \{1,...,N(\varepsilon)\}$. In view of the same property, the Young, H\"older and Sobolev inequalities in $D^n_1$ provide
\begin{equation} \label{catta1}
\begin{aligned}
\| a^{2}_{\varepsilon} F - a_{\varepsilon} (U_{\varepsilon}\cdot\nabla)U_{\varepsilon} \|^{3/2}_{L^{3/2}(D^n_1)} & \leq C \left( a^{3}_{\varepsilon} \|  F  \|^{3/2}_{L^{3/2}(D^n_1)} + a^{3/2}_{\varepsilon} \| (U_{\varepsilon}\cdot\nabla)U_{\varepsilon} \|^{3/2}_{L^{3/2}(D^n_1)} \right) \\[5pt]
& \leq C \left( a^{3}_{\varepsilon} \|  F  \|^{3/2}_{L^{3/2}(D^n_1)} + a^{3/2}_{\varepsilon} \| \nabla U_{\varepsilon} \|^{3/2}_{L^{2}(D^n_1)} \| U_{\varepsilon} \|^{3/2}_{L^{6}(D^n_1)} \right) \\[5pt]
& \leq C \left( a^{3}_{\varepsilon} \|  F  \|^{3/2}_{L^{3/2}(D^n_1)} + a^{3/2}_{\varepsilon} \| \nabla U_{\varepsilon} \|^{3}_{L^{2}(D^n_1)} \right) \, .
\end{aligned}
\end{equation}
Since $a_{\varepsilon} \ll \varepsilon$ for every $\varepsilon \in I_{*}$, we clearly have
$$
D^n_1 \subset \widetilde{D^n_1} \doteq D \left( (0,0), \dfrac{\delta_{1} \varepsilon}{a_{\varepsilon}} \right) \setminus \overline{K_{n}} \, ,
$$
so that successive applications of the change of variables \eqref{cv} ensure that
\begin{equation} \label{catta2}
\begin{aligned}
& \|  F  \|^{3/2}_{L^{3/2}(D^n_1)} \leq \dfrac{1}{a^{2}_{\varepsilon}} \|  f  \|^{3/2}_{L^{3/2}(D \left(\xi^{\varepsilon}_{n}, \delta_{1} \varepsilon \right) \setminus \overline{K_{n}^{\varepsilon}})} \, , \qquad \|  \nabla U_{\varepsilon}  \|_{L^{2}(D^n_1)} \leq \|  \nabla u_{\varepsilon}  \|_{L^{2}(D \left(\xi^{\varepsilon}_{n}, \delta_{1} \varepsilon \right) \setminus \overline{K_{n}^{\varepsilon}})} \, , \\[6pt]
& \|  \nabla U_{\varepsilon}   \|^{3/2}_{L^{3/2}(D^n_1)} \leq \dfrac{1}{\sqrt{a_{\varepsilon}}} \|  \nabla u_{\varepsilon}   \|^{3/2}_{L^{3/2}(D \left(\xi^{\varepsilon}_{n}, \delta_{1} \varepsilon \right) \setminus \overline{K_{n}^{\varepsilon}})} \, , \qquad \|  P_{\varepsilon}  \|^{3/2}_{L^{3/2}(D^n_1)} \leq \dfrac{1}{a^{2}_{\varepsilon}} \|  p_{\varepsilon}  \|^{3/2}_{L^{3/2}(D \left(\xi^{\varepsilon}_{n}, \delta_{1} \varepsilon \right) \setminus \overline{K_{n}^{\varepsilon}})} \, .
\end{aligned}
\end{equation}
Inserting \eqref{catta1}-\eqref{catta2} into \eqref{regularity1} gives us
\begin{equation} \label{regularity2}
\begin{aligned}	
	& \| \nabla p_{\varepsilon} \|^{3/2}_{L^{3/2}(D \left(\xi^{\varepsilon}_{n}, \delta_{0} a_{\varepsilon} \right) \setminus \overline{K_{n}^{\varepsilon}})} = \sqrt{a_{\varepsilon}} \, \| \nabla P_{\varepsilon} \|^{3/2}_{L^{3/2}(D^n_0)} \leq C \Bigg(  a^{3/2}_{\varepsilon} \|  f  \|^{3/2}_{L^{3/2}(D \left(\xi^{\varepsilon}_{n}, \delta_{1} \varepsilon  \right) \setminus \overline{K_{n}^{\varepsilon}})}  \\[6pt]
	& \hspace{12mm} + a^{2}_{\varepsilon} \| \nabla u_{\varepsilon} \|^{3}_{L^{2}(D \left(\xi^{\varepsilon}_{n}, \delta_{1} \varepsilon  \right) \setminus \overline{K_{n}^{\varepsilon}})} + \|  \nabla u_{\varepsilon}  \|^{3/2}_{L^{3/2}(D \left(\xi^{\varepsilon}_{n}, \delta_{1} \varepsilon  \right) \setminus \overline{K_{n}^{\varepsilon}})} + \dfrac{1}{a^{3/2}_{\varepsilon}} \| p_{\varepsilon} \|^{3/2}_{L^{3/2}(D \left(\xi^{\varepsilon}_{n}, \delta_{1} \varepsilon  \right) \setminus \overline{K_{n}^{\varepsilon}})} \Bigg) \, .
\end{aligned}	
\end{equation}
Following the proof of \cite[Theorem 2.1]{lu2021uniform}, we decompose the perforated domain as
$$
\Omega_{\varepsilon} = \left[ \Omega_{\varepsilon} \setminus \left( \bigcup_{n=1}^{N(\varepsilon)} D \left(\xi^{\varepsilon}_{n}, \delta_{0} a_{\varepsilon} \right) \right) \right] \cup \bigcup_{n=1}^{N(\varepsilon)} \left( D \left(\xi^{\varepsilon}_{n}, \delta_{0} a_{\varepsilon} \right) \setminus \overline{K_{n}^{\varepsilon}} \right) \doteq \widetilde{\Omega_{\varepsilon}} \cup \bigcup_{n=1}^{N(\varepsilon)} \left( D \left(\xi^{\varepsilon}_{n}, \delta_{0} a_{\varepsilon} \right) \setminus \overline{K_{n}^{\varepsilon}} \right) \, .
$$
On one hand, since the holes are mutually disjoint (see \eqref{perforation}$_3$), from \eqref{number0}-\eqref{regularity2} we directly obtain
\begin{equation} \label{masmoudi1}
\begin{aligned}
& \| \nabla p_{\varepsilon} \|^{3/2}_{L^{3/2}(\bigcup_{n=1}^{N(\varepsilon)} \left( D \left(\xi^{\varepsilon}_{n}, \delta_{0} a_{\varepsilon} \right) \setminus \overline{K_{n}^{\varepsilon}} \right))} = \sum_{n=1}^{N(\varepsilon)} \| \nabla p_{\varepsilon} \|^{3/2}_{L^{3/2}(D \left(\xi^{\varepsilon}_{n}, \delta_{0} a_{\varepsilon} \right) \setminus \overline{K_{n}^{\varepsilon}})} \\[6pt]
& \hspace{-4mm} \leq C \Bigg(  a^{3/2}_{\varepsilon} \|  f  \|^{3/2}_{L^{3/2}(\Omega)} + \dfrac{a^{2}_{\varepsilon}}{\varepsilon^{2}} \| \nabla u_{\varepsilon} \|^{3}_{L^{2}(\Omega_{\varepsilon})} + \|  \nabla u_{\varepsilon}  \|^{3/2}_{L^{3/2}(\Omega_{\varepsilon})} + \dfrac{1}{a^{3/2}_{\varepsilon}} \| p_{\varepsilon} \|^{3/2}_{L^{3/2}(\Omega_{\varepsilon})} \Bigg) \\[6pt]
& \hspace{-4mm} \leq \dfrac{C}{a^{3/2}_{\varepsilon}} \left(  1 + \|  f  \|^{2}_{L^{2}(\Omega)} + \|  \nabla u_{\varepsilon}  \|^{4}_{L^{2}(\Omega_{\varepsilon})} + \| p_{\varepsilon} \|^{2}_{L^{2}(\Omega_{\varepsilon})} \right)^{3/2} \qquad \forall \varepsilon \in I_{*} \, .
\end{aligned}	
\end{equation}
On the other hand, since the points of the region $\widetilde{\Omega_{\varepsilon}}$ are sufficiently far away from the holes $K_{\varepsilon}$, we can argue as in \cite[Theorem 3.3]{gazspe} to obtain the bound
\begin{equation} \label{masmoudi2}
	\| \nabla p_{\varepsilon} \|_{L^{3/2}(\widetilde{\Omega_{\varepsilon}})} \leq \dfrac{C}{a_{\varepsilon}} \left(1 + \| f \|^{2}_{L^{2}(\Omega)} + \| v_{*} \|^{4}_{H^{3/2}(\partial \Omega)} \right) \qquad \forall \varepsilon \in I_{*} \, .
\end{equation}
Adding \eqref{masmoudi1}-\eqref{masmoudi2} gives us \eqref{masmoudi0}. It is therefore left open the possibility of recovering the results of Theorems \ref{epslevel}-\ref{effectiveq1} without the smallness assumption \eqref{neustupa}.

\par\medskip\noindent
{\bf Data availability statement.} Data sharing not applicable to this article as no datasets were generated or analyzed during the current study.
\par\smallskip
\noindent
{\bf Conflict of interest statement}.  The authors declare that there is no conflict of interest.

\phantomsection
\addcontentsline{toc}{section}{References}
\bibliographystyle{abbrv}
\bibliography{references}
\vspace{5mm}

\noindent
\hspace{0.1mm}
\begin{minipage}{140mm}
	\textbf{Clara Patriarca}\\
	Département de Mathématique\\
	Université Libre de Bruxelles\\
	 Boulevard du Triomphe 155\\
	1050 Brussels - Belgium\\
	E-mail: clara.patriarca@ulb.be
\vspace{0.8cm}		
\end{minipage}
\newline
\vspace{0.8cm}
\noindent
\begin{minipage}{100mm}
\textbf{Gianmarco Sperone}\\
Dipartimento di Matematica\\
Politecnico di Milano\\
Piazza Leonardo da Vinci 32\\
20133 Milan - Italy\\
E-mail: gianmarcosilvio.sperone@polimi.it\\
\textit{currently at}\\
Facultad de Matemáticas\\
Pontificia Universidad Católica de Chile\\
Avenida Vicuña Mackenna 4860\\
7820436 Santiago - Chile\\
E-mail: gianmarco.sperone@uc.cl
\end{minipage}
\end{document}